\documentclass[11pt]{article}

\usepackage{tikz}
\usepackage{enumerate}
\usepackage{booktabs}
\usepackage{bookmark}
\usepackage{geometry}

\usepackage{mathtools}
\usepackage{mathrsfs}
\usepackage{dsfont}
\usepackage{amssymb}
\usepackage{amsthm}

\usepackage{bm} 
\usepackage{graphics}
\makeatletter
\DeclareRobustCommand*\uell{\mathpalette\@uell\relax}
\newcommand*\@uell[2]{
    \setbox0=\hbox{$#1\ell$}
    \setbox1=\hbox{\rotatebox{10}{$#1\ell$}}
    \dimen0=\wd0 \advance\dimen0 by -\wd1 \divide\dimen0 by 2
    \mathord{\lower 0.1ex
    \hbox{\kern\dimen0\unhbox1\kern\dimen0}}
}
\makeatother

\makeatletter
\newcommand*{\@restrictionaux}[2]{{#1\,\smash{\vrule height 0.7\ht1 depth 1.2\dp1}}_{\,#2}} \newcommand*{\restr}[2]{\mathchoice
  {\setbox1\hbox{${\displaystyle #1}_{\scriptstyle #2}$} \@restrictionaux{#1}{#2}}
  {\setbox1\hbox{${\textstyle #1}_{\scriptstyle #2}$} \@restrictionaux{#1}{#2}}
  {\setbox1\hbox{${\scriptstyle #1}_{\scriptscriptstyle #2}$} \@restrictionaux{#1}{#2}}
  {\setbox1\hbox{${\scriptscriptstyle #1}_{\scriptscriptstyle #2}$} \@restrictionaux{#1}{#2}}
}
\makeatother

\newcommand{\mysetminusD}{\hbox{\tikz{\draw[line width=0.6pt,line cap=round] (3pt,0) -- (0,6pt);}}}
\newcommand{\mysetminusT}{\mysetminusD}
\newcommand{\mysetminusS}{\hbox{\tikz{\draw[line width=0.45pt,line cap=round] (2pt,0) -- (0,4pt);}}}
\newcommand{\mysetminusSS}{\hbox{\tikz{\draw[line width=0.4pt,line cap=round] (1.5pt,0) -- (0,3pt);}}}
\newcommand{\mysetminus}{\mathbin{\mathchoice{\mysetminusD}{\mysetminusT}{\mysetminusS}{\mysetminusSS}}}
\renewcommand{\setminus}{\mysetminus}

\renewcommand{\phi}{\varphi}
\renewcommand{\epsilon}{\varepsilon}

\renewcommand{\P}{\mathds{P}}
\newcommand{\E}{\mathds{E}}
\newcommand{\N}{\mathds{N}}
\newcommand{\Q}{\mathds{Q}}
\newcommand{\R}{\mathds{R}}
\newcommand{\M}{\mathrm{M}}

\newcommand{\goesto}[2][]{\xrightarrow[\:#2\:]{\:#1\:}}
\newcommand{\indic}{\mathbf{1}}
\newcommand{\diff}{\mathrm{d}}
\DeclarePairedDelimiter{\abs}{|}{|}
\DeclarePairedDelimiter{\norm}{\Vert}{\Vert}
\DeclarePairedDelimiter{\Angle}{\langle}{\rangle}
\renewcommand{\angle}{\Angle}
\DeclarePairedDelimiter{\floor}{\lfloor}{\rfloor}
\DeclareMathOperator{\Card}{Card}
\DeclareMathOperator{\Leb}{Leb}
\DeclareMathOperator{\supp}{supp}
\DeclareMathOperator{\h}{ht}

\newtheorem{theorem}{Theorem}[section]
\newtheorem{proposition}[theorem]{Proposition}
\newtheorem{lemma}[theorem]{Lemma}
\newtheorem{corollary}[theorem]{Corollary}

\theoremstyle{definition}

\newtheorem{assumption}{Assumption}
\newtheorem{remark}{Remark}

\usepackage{xcolor}
\definecolor{myteal}{HTML}{004b6f}
\definecolor{myblack}{HTML}{000000}
\definecolor{citegrey}{HTML}{888888}
\usepackage{hyperref}
\hypersetup{
    colorlinks=true,
    linkcolor=myteal,
    citecolor=citegrey,
    urlcolor=myteal,
}

\usepackage{authblk}

\title{
    A moment approach for the convergence of spatial branching processes
    to the Continuum Random Tree
}
\author[1]{Félix Foutel-Rodier}
\date{}
\affil[1]{Department of Statistics, University of Oxford}

\begin{document}

\maketitle

\begin{abstract}
    We consider a general class of branching processes in discrete time,
    where particles have types belonging to a Polish space and reproduce
    independently according to their type. If the process is critical and
    the mean distribution of types converges for large times, we prove
    that the tree structure of the process converges to the Brownian
    Continuum Random Tree, under a moment assumption.
    We provide a general approach to prove similar invariance principles
    for branching processes, which relies on deducing the convergence of the
    genealogy from computing its moments. These are obtained using a new
    many-to-few formula, which provides an expression for the moments of order $k$
    of a branching process in terms of a Markov chain indexed by a uniform tree
    with $k$ leaves.
\end{abstract}

\section{Introduction}

\subsection{Context and motivation}

This article is concerned with scaling limits of the tree structure of
spatial branching processes. In discrete time, a spatial branching
process is a particle system where each particle has a type in some
Polish type space $E$, which we think of as a set of spatial
locations. At each generation a particle at $x \in E$ is replaced by its
offspring, encoded as a point process
\begin{equation} \label{eq:reproductionPP}
    \Xi_x = \sum_{i=1}^{\abs{\Xi_x}} \delta_{\xi_{x,i}},
\end{equation}
with the interpretation that the parent has $\abs{\Xi_x}$ children,
located at $(\xi_{x,i})_{i \le \abs{\Xi_x}}$. The next generation is
obtained by carrying out this step independently for each particle.
This formalism encompasses many well-studied branching models such as
multi-type branching processes, branching random walks, or any branching
diffusion observed at discrete time points 
\cite{Athreya1972, shi2015branching, pinsky1995positive}. 
Envisioning the particles as vertices and connecting each vertex to its
parent, this procedure constructs a random tree $T$ augmented with a
collection of random variables $(X_v)_{v \in T}$ giving the types. 
Our objective is to study scaling limits of the marked tree $(T,
(X_v)_{v \in T})$.

Such scaling limits have already been investigated extensively in some
special cases. In the absence of types, that is when $E$ is made of a
single element, we recover Bienaymé--Galton--Watson processes. Scaling
limits of their tree structure have received considerable attention,
starting with seminal work of Aldous on the Brownian Continuum Random
Tree (CRT) \cite{aldous1993continuum} and subsequently leading to a full
classification in the near-critical case and the introduction of Lévy
trees \cite{legall1999spatial, duquesne2002random}. One can extend these
results to a special family of spatial branching processes by assuming
that the number of offspring $\abs{\Xi_x}$ in \eqref{eq:reproductionPP}
does not depend on $x$. In this case, the tree structure $T$ remains
distributed as a Bienaymé--Galton--Watson tree and the spatial locations
$(X_v)_{v \in T}$ can be obtained by running a Markov chain indexed by
the vertices of that tree. It is then possible to describe the joint
scaling limit of the genealogy and particle locations as a Markov process
indexed by a continuous tree, which is tightly connected to the notion of
snake \cite{LeGall1991, LeGall1993_2, LeGall1998} and corresponds to the
genealogy of the so-called Dawson--Watanabe superprocesses
\cite{dawson1993measurevalued, perkins1999, etheridge2000introduction}.
An important example of such processes is the branching Brownian motion
and its scaling, the super-Brownian motion.

However, for general offspring distributions $T$ fails to be
a Bienaymé--Galton--Watson tree and known results become much scarcer,
with a few notable exceptions \cite{miermont2008invariance, lambert2010,
berestycki2013, powell2019invariance}. In this case, the tree structure
$T$ and the particle locations $(X_v)_{v \in T}$ become strongly
dependent: vertices carrying a successful type are typically expected to
have a larger degree and more descendants. Several techniques based on
encoding trees as paths become difficult to apply because of this
dependence, making the analysis of $T$ challenging. Beyond the
technical challenge, letting the number of descendants depend on the
particles location leads to a much richer class of random trees $T$ and
opens the possibility of observing new interesting behaviour emerging
from the addition of a spatial structure. The references above already
provide interesting examples of \emph{binary} branching diffusions where
the presence of highly successful types drives the system to the
universality class of  $\alpha$-stable branching processes
\cite{berestycki2013, Tourniaire2024, foutel2024convergence},
with the emergence of infinite branch points in the limit. This is in
sharp contrast with Bienaymé--Galton--Watson trees where, in the absence
of type, $\alpha$-stable branching processes can only arise if the
offspring distribution has heavy tail. A general theory for these spatial
branching processes is still missing, which seems to be an interesting
avenue and motivates the introduction of new tools to study their tree
structures.

Quite recently, starting with the work of Harris and Roberts
\cite{Harris2017}, some new techniques have been developed 
under the name of \emph{many-to-few formulas} to study the genealogical
properties of general spatial branching processes. They rely on
expressing the law of the subtree spanned by $k$ typical particles in the
population using an appropriate random change of measure, in the spirit
of the celebrated spinal decomposition results for branching processes
\cite{chauvin1991growing, lyons1995conceptual, biggins2004measure}. These
techniques seem very promising and have already been applied successfully
in various contexts \cite{harris2022coalescent, Johnston2019,
harris2024universality, foutel2023convergence, boenkost2022genealogy,
schertzer2023spectral, foutel2024convergence}. 
We can now express the main two goals of the present work.
\begin{itemize}
    \item Provide a general approach to study scaling limits of the marked
        tree $(T, (X_v)_{v \in T})$ that relies on a many-to-few formula
        and a method of moments.
    \item Prove an invariance principle to the Brownian CRT for a large
        class of critical spatial branching processes with an ergodic
        mean behaviour.
\end{itemize}
We describe the invariance principle and the assumptions under which it
holds in the next section. Our moment approach extends that proposed in
\cite{foutel2023convergence}, which is restricted to studying the
\emph{reduced tree} spanned by a single generation. It relies on an
adaptation of the method of moments to infinite measures developed in
\cite{foutel2024vague} and on a new many-to-few formula, which is another
contribution of our work. It is described informally in
Section~\ref{SS:proofIdea} and made precise in Section~\ref{S:moments}.
We also refer to the recent work of \cite{Harris2024} for a different
expression in the same direction.

Before ending this section let us mention that, in addition to being a
natural probabilistic problem, extending our understanding of the tree
structure of spatial branching processes beyond the class of tree-indexed
Markov chains is strongly motivated by applications. Types can be used to
model a wide variety of features of interest, ranging from spatial
locations, traits, and age in population biology \cite{Baird2003,
calvez2022, jagers1984, etheridge2022genealogies}, clinical status in
epidemiology \cite{britton2019estimation, foutelrodier2022}, to positions
and velocities in neutron transport \cite{horton2023stochastic}. Modeling
most real-world phenomena requires to have the offspring law depend on
the location, and studying the resulting tree structure often leads to a
deeper understanding of the dynamics of the process. This is especially
true in population genetics, where genealogies are key to analysing the
patterns of genetic diversity. 
In addition to that, trees and marked trees are central to the study of many
combinatorial objects \cite{Chassaing2004, Marckert2007, addarioberry2012}, 
an important example being planar maps, and progress in our understanding
of spatial branching processes might also be relevant in that direction.

\subsection{Invariance principle for critical spatial branching processes} 

The invariance principle that we derive holds for a class of spatial
branching processes such that the average particle location reaches an
equilibrium. Our assumptions are almost directly
borrowed from \cite{gonzalez2022asymptotic, harris2022yaglom}, which
themselves trace back to early work of Asmussen and Hering on general
branching processes \cite{hering1971, Asmussen1983, Asmussen1976}. We
will use the notation $\mu[f]$ or $\angle{\mu,f}$ for the integral of $f$
against some measure $\mu$, and let $\P_x$ be the law of the process
started from a single ancestor at $x \in E$.

\begin{assumption} \label{ass:criticality}
    The following three points hold.
    \begin{enumerate}[(i)]
    \item    
    There exists a continuous bounded harmonic function $h \colon E \to
    (0, \infty)$ and a stationary measure $\pi$ on $E$ for the mean
    semi-group of the process, namely
    \[
        \forall x \in E,\quad h(x) = \E\big[\angle{\Xi_x, h} \big],
        \qquad
        \pi(\diff x) = \int_E \E[ \Xi_y(\diff x) ] \pi(\diff y).
    \]
    They are renormalised so that $\pi$ is a probability measure and
    $\angle{\pi, h} = 1$.
    \item For any continuous bounded $f \colon E \to \R$,
    \begin{equation} \label{eq:convMean}
        \adjustlimits \lim_{n \to \infty} \sup_{x \in E} \: 
        \abs[\Big]{\E_x\Big[ \sum_{v \in T, \abs{v} = n} f(X_v) \Big] - h(x) \angle{\pi, f} } = 0.
    \end{equation}
    \item For any $k \ge 1$, 
        \begin{equation} \label{eq:momentBound}
            \sup_{x \in E} \E\big[ \abs{\Xi_x}^k \big] < \infty
        \end{equation}
        and the map $x \mapsto \E\big[ \abs{\Xi_x}^2 \big]$ is continuous.
    \end{enumerate}
\end{assumption}

Point (i) provides the existence of a pair $(h, \pi)$ of right
eigenfunction and left eigenmeasure for the mean semi-group corresponding
to the eigenvalue $1$. Point (ii) is a stability condition for this
semi-group, requiring that it has a Perron--Frobenius-like behaviour as
$n \to \infty$. These two conditions serve as a notion of criticality for
our spatial branching processes. The measure $\pi$ corresponds to the
asymptotic distribution of particles, and the harmonic function $h$
quantifies the ``success'' of a type: a particle of type $x$ will have
asymptotically on average $h(x)$ descendants. Point (iii) will be
required for applying our method of moments. We introduce
\[
    \Sigma^2 \coloneqq \int_E \E\bigg[ \sum_{\substack{i,j=1\\ i \ne j}}^{\abs{\Xi_x}}
    h(\xi_{x,i}) h(\xi_{x,j}) \bigg] \, \pi(\diff x),
\]
which is the key parameter driving the dynamics of the process. It acts
as a spatial notion of variance and is clearly finite under
Assumption~\ref{ass:criticality}. 

Let us illustrate rapidly these conditions with a few examples. If $E$ is
finite as in \cite{miermont2008invariance}, the Perron--Frobenius theorem
readily ensures that (i) and (ii) are fulfilled whenever the process is
critical, irreducible, and aperiodic. A second example is obtained by
observing a branching diffusion killed upon reaching the boundary of a
domain $E \subseteq \R^d$ at discrete time steps. (Continuous or discrete
time makes no difference in our approach, but we chose the latter because
it is  notationally more convenient for manipulating trees.) If the
domain $E$ is bounded, spectral theory for elliptic operators entails
that (i) and (ii) hold if the process is critical
\cite{pinsky1995positive, powell2019invariance}. These points also hold
for some branching diffusions on unbounded domains if the potential
decays sufficiently fast at infinity \cite{collet2024branching}. Two last
examples that we want to mention are critical growth--fragmentation
processes \cite{bertoin2018probabilistic} and stochastic neutron
transport \cite{horton2023stochastic}, the latter being one the original
motivation of \cite{harris2022yaglom}. Finally, see \cite{bansaye2022non,
delMoral2023stability} for general results allowing one to establish
\eqref{eq:convMean}.

Assumption~\ref{ass:criticality} is sufficient to derive a first
version of our invariance principle. In order to discuss scaling limits
of trees, we need to view our branching process as a
random \emph{marked metric measure space}. The appropriate state space
and topologies are introduced formally in Section~\ref{S:mmmSpaces}.
Recall that a spatial branching process defines a random tree $T$
decorated with a collection of types $(X_v)_{v \in T}$, and let $T$ be
rooted at the initial particle $\rho$. The tree $T$ is naturally endowed
with its graph distance $d$ and with a measure $\mu$ on $T \times E$
defined as 
\[
    \mu = \sum_{v \in T} \delta_{(v, X_v)}.
\]
This measure encodes both the number of individuals and their types.
In the terminology of \cite{Depperschmidt2011}, $\mathcal{T} = (T, d,
\rho, \mu)$ is a (pointed) marked metric measure space representing
the genealogy and types of the population. We rescale edge lengths by $n$
and particle masses by $n^2$ to define 
\begin{equation} \label{eq:rescalingMmmSpace}
    \bar{\mathcal{T}}_n \coloneqq \big( T, \tfrac{d}{n}, \rho,
    \tfrac{\mu}{n^2} \big),
\end{equation}
and let $\mathscr{L}_x(\bar{\mathcal{T}}_n)$ denote the law of this
rescaled metric space, when starting from a single particle at $x$.
Finally, we let $\mathcal{T}_{\mathrm{b}, \pi} = (T_\mathrm{b}, d_{\mathrm{b}},
\rho_{\mathrm{b}}, \mu_{\mathrm{b}} \otimes \pi)$ denote a free Brownian CRT
with variance $\Sigma^2$ and marks $\pi$, which is the tree encoded by a
Brownian excursion with variance $4/\Sigma^2$ whose leaves are marked
independently according to $\pi$. This object is constructed precisely in
Section~\ref{SS:BrownianCRT}. We let $\mathscr{L}(\mathcal{T}_{\mathrm{b},\pi})$ 
be its ``distribution,'' which is the law of this random variable under
the infinite excursion measure of Brownian motion.

Even if it is not our main focus here, it comes at no cost with our
approach to derive the scaling limit of a second tree: the \emph{reduced
tree} spanned by a given generation. Fix $n \ge 0$ and let us denote by
$U_n$ the $n$-th generation of the process defined as 
\[
    U_n \coloneqq \{ v \in T : \abs{v} = n \}.
\]
Let $d_n$ be the restriction of the graph distance
$d$ to $U_n$ and $\tilde{\mu}_n$ be the empirical measure of labels and
types at generation $n$, namely
\[
    \tilde{\mu}_n = \sum_{v \in U_n} \delta_{(v, X_v)}.
\]
The marked metric measure space $(U_n, d_n, \rho, \tilde{\mu}_n)$ is an
\emph{ultrametric space} that encodes the genealogy and types of the
$n$-th generation of the spatial branching process. We rescale it as
above to define
\begin{equation} \label{eq:defUMS}
    \bar{\mathcal{U}}_n \coloneqq \big( U_n, \tfrac{d_n}{n}, 
    \rho, \tfrac{\tilde{\mu}_n}{n} \big).
\end{equation}
Note the change in mass rescaling compared to \eqref{eq:rescalingMmmSpace}.
Let $\mathcal{U}_{\mathrm{b}, \pi} = (U_{\mathrm{b},\pi},
d_{\mathrm{b},\pi}, \rho_{\mathrm{b},\pi}, \mu_{\mathrm{b}}\otimes \pi)$ be the
Brownian coalescent point process (CPP) \cite{Popovic2004, Lambert2013}
with variance $\Sigma^2$ and marks $\pi$. This ultrametric space
corresponds to the sphere of radius $1$ of a Brownian CRT conditioned to
have height at least $1$, see Section~\ref{SS:BrownianCRT} for a precise
definition. Finally, let us denote by $Z_n = \abs{U_n}$ the number of
particles at generation $n \ge 0$.

\newpage

\begin{theorem} \label{thm:GromovVague}
    Suppose that Assumption~\ref{ass:criticality} is fulfilled.
    \begin{enumerate}[\rm (i)]
        \item The convergence
            \begin{equation*} 
                \forall x \in E,\quad n \mathscr{L}_x\big( \bar{\mathcal{T}}_n \big)
                \goesto{n \to \infty}
                h(x) \mathscr{L}(\mathcal{T}_{\mathrm{b},\pi})
            \end{equation*}
            holds vaguely in the marked Gromov-vague topology, where
            $\mathcal{T}_{\mathrm{b}, \pi}$ is a free Brownian CRT with variance
            $\Sigma^2$ and independent marks $\pi$.
        \item There exists $(\eta_n)_{n \ge 1}$ with $\eta_n \to 0$
            such that
            \begin{equation} \label{eq:weakKolmogorov}
                \forall x \in E,\quad 
                \lim_{n \to \infty} n \P_x( Z_n \ge \eta_n n ) \goesto{n \to \infty}
                \frac{2h(x)}{\Sigma^2}.
            \end{equation}
            Moreover, conditional on $\{ Z_n > n\eta_n \}$ and starting
            from one particle at $x \in E$,
            \[
                \bar{\mathcal{U}}_n
                \goesto{n \to \infty} 
                \mathcal{U}_{\mathrm{b}, \pi}
            \]
            in distribution for the marked Gromov-weak topology, where
            $\mathcal{U}_{\mathrm{b}, \pi}$ is a Brownian CPP with variance
            $\Sigma^2$ and marks $\pi$.
    \end{enumerate}
\end{theorem}

\begin{remark}
    We conjecture that the assumptions of this result can be relaxed as
    long as $\Sigma^2 < \infty$. In particular $h$ need not be bounded
    and $\abs{\Xi_x}$ need only have finite moments of order $k=2$. Our
    approach could lead to the latter improvement with some additional
    effort by truncating the process as in \cite{boenkost2022genealogy},
    but we restrain from doing so to keep the proofs more transparent. 
    Conversely, we believe that relaxing the former assumption would lead
    to new behaviours. If $h$ is unbounded, it is possible that
    $\Sigma^2=\infty$ even if $\abs{\Xi_x}$ has finite variance. Several
    examples have shown that, in this case, $\alpha$-stable genealogies
    can arise due to rare excursions of particles to infinity, generating
    a large offspring in a small number of generations \cite{berestycki2013,
    Tourniaire2024, foutel2024convergence}. We believe this to
    be a generic behaviour of spatial branching processes and that
    relaxing our assumption to include unbounded $h$ is an interesting
    avenue to investigate.
\end{remark}

The Gromov-vague convergence obtained above is equivalent to the
Gromov-weak convergence (also known as Gromov--Prohorov convergence) of
the tree below a height of order $n$. This result already covers many
interesting applications, including deriving the limit in law of the
subtree spanned by $k$ individuals sampled uniformly from the population
before or at a fixed generation, in the spirit of coalescent theory in
population genetics \cite{kingman1982coalescent, Berestycki2009}.
It also provides a weak estimate on the survival probability
\eqref{eq:weakKolmogorov}, which gives the probability that the
population reaches a ``macroscopic'' size of order $n$.

However, the topologies used above are too weak for other applications, as
many natural geometric quantities (such as the diameter and the total
height) are not continuous functionals of the Gromov-weak topology. This
provides motivation for reinforcing this result to the stronger
Gromov--Hausdorff--Prohorov topology. In the context of spatial branching
processes, we will show that reinforcing the convergence amounts to
sharpening the estimate \eqref{eq:weakKolmogorov} on the survival
probability. This leads us to formulate a second assumption on our
branching process, which should be thought of as reinforcing
\eqref{eq:weakKolmogorov}.

\begin{assumption} \label{ass:kolmogorov}
    Suppose that 
    \begin{equation} \label{eq:kolmogorov}
        \adjustlimits\lim_{n \to \infty} \sup_{x \in E} \: 
        \abs[\Big]{ n\P_x(Z_n > 0) - \frac{2h(x)}{\Sigma^2} }
        = 0.
    \end{equation}
\end{assumption}

\begin{remark}
    In a slightly different setting, \cite{harris2022yaglom, horton2024stability} 
    have shown that \eqref{eq:kolmogorov} follows from Assumption~\ref{ass:criticality} 
    provided that the process goes extinct almost surely and verifies a
    mild irreducibility condition. In our notation, the latter condition
    is that there are some $K, M > 0$ such that
    \[
        \forall f \colon E\to \R_+,\quad 
        \int_E \E_x\Big[ 
            \indic_{\{ \abs{\Xi} \le M\}} \sum_{i, j = 1, i\ne j}^{\abs{\Xi}} 
            f(\xi_i)f(\xi_j) \Big] \pi(\diff x)
        \ge K \angle{\pi, f}.
    \]
    We believe that the proof of \cite{harris2022yaglom, horton2024stability} 
    could be adapted to the present setting. Since this is not our focus
    here, we prefer to leave \eqref{eq:kolmogorov} as an additional
    condition to verify rather than reproducing the arguments in
    \cite{horton2024stability}.
\end{remark}

For Bienaymé--Galton--Watson processes, the asymptotic
\eqref{eq:kolmogorov} is known as Kolmogorov's estimate. More generally,
it is expected to hold for critical branching processes with a finite
``variance.'' For instance, \eqref{eq:kolmogorov} is known to hold in
most of the examples mentioned above. With Kolmogorov's estimate at hand,
we can reinforce Theorem~\ref{thm:GromovVague} to the following result.

\begin{theorem} \label{thm:GHP}
    Suppose that Assumption~\ref{ass:criticality} and
    Assumption~\ref{ass:kolmogorov} are fulfilled. 
    \begin{enumerate}[\rm (i)]
        \item The convergence 
            \begin{equation*} 
                \forall x \in E,\quad n \mathscr{L}_x\big( \bar{\mathcal{T}}_n \big)
                \goesto{n \to \infty}
                h(x) \mathscr{L}(\mathcal{T}_{\mathrm{b},\pi})
            \end{equation*}
            holds vaguely in the marked Gromov--Hausdorff--Prohorov topology.
        \item Conditional on $\{ Z_n > 0 \}$ and starting from one particle at $x \in E$, 
            \[
                \bar{\mathcal{U}}_n
                \goesto{n \to \infty} 
                \mathcal{U}_{\mathrm{b}, \pi}
            \]
            in distribution for the marked Gromov--Hausdorff--Prohorov
            topology.
    \end{enumerate}
\end{theorem}

Point~(i) of the above result might remain daunting as it involves
infinite measures and vague convergence. As a direct corollary, we can
deduce two more standard types of convergence to the Brownian CRT. The
first one is the convergence of the branching process conditioned to
survive for a long time $tn$ to the Brownian CRT conditioned to have
height larger than $t$. The second one is the convergence of the forest
obtained by starting the population from $n$ ancestors to the forest
encoded by a reflected Brownian motion. Let us introduce some notation
for the latter case. Suppose that the population starts from a set of
$z_0$ ancestors with positions $(x_i)_{i \le z_0}$, which we record as a
measure $\mathbf{z}_0 = \delta_{x_1}+\dots+\delta_{x_{z_0}}$.
Each ancestor grows an independent tree, and we denote by 
$(\bar{\mathcal{T}}_{n,i})_{i \le z_0}$ the collection of metric
measure spaces obtained by applying the rescaling \eqref{eq:rescalingMmmSpace} 
to each tree and re-ordering them in decreasing order of total mass,
\[
    \abs{ \bar{\mathcal{T}}_{n,1} } 
    \ge \abs{ \bar{\mathcal{T}}_{n,2} } 
    \ge \dots,
\]
breaking ties arbitrarily. We let $\h(\mathcal{T})$ denote the height 
of a tree $\mathcal{T}$, see \eqref{eq:heightDef}.

\begin{corollary} \label{cor:GHPconsequence}
    Suppose that Assumption~\ref{ass:criticality} and
    Assumption~\ref{ass:kolmogorov} are fulfilled.
    \begin{enumerate}[\rm (i)]
        \item For any $t > 0$ and $x \in E$, 
            \[
                \mathscr{L}_x(\bar{\mathcal{T}}_n \mid Z_{tn}
                > 0 )
                \goesto{n \to \infty}
                \mathscr{L}\big( \mathcal{T}_{\mathrm{b},\pi} \mid
                \h(\mathcal{T}_{\mathrm{b},\pi}) \ge t \big)
            \]
            in distribution for the marked Gromov--Hausdorff--Prohorov topology.
        \item Suppose that the population starts from an initial
            configuration of particles $\mathbf{z}_{0,n}$ 
            such that $\frac{1}{n} \mathbf{z}_{0,n} \to \nu_0$ weakly as
            $n \to \infty$ for some finite measure $\nu_0$ on $E$. Let
            $(\mathcal{T}_{\mathrm{b},i})_{i \ge 1}$ be the atoms of a
            Poisson point process with intensity measure $\angle{\nu_0,
            h} \mathscr{L}(\mathcal{T}_{\mathrm{b},\pi})$, in decreasing
            order of their total mass. Then
            \[
                (\bar{\mathcal{T}}_{n,1}, \bar{\mathcal{T}}_{n,2}, \dots)
                \goesto{n \to \infty}
                (\mathcal{T}_{\mathrm{b},1}, \mathcal{T}_{\mathrm{b},2}, \dots)
            \]
            in distribution in the finite-dimensional sense, where each
            coordinate is endowed with the marked Gromov--Hausdorff--Prohorov
            topology.
    \end{enumerate}
\end{corollary}

\subsection{Main proof ideas: moments and many-to-few formula} 
\label{SS:proofIdea}

Our proof of the invariance principle relies on two main ingredients.
First, we will deduce the convergence in the Gromov-vague topology from a
method of moments. The moments of order $k \ge 1$ of the random marked tree 
$\mathcal{T}$ are quantities of the form 
\begin{equation} \label{eq:momentIntro}
    \E_x[\Phi(\mathcal{T})]
    \coloneqq \E_x\Big[ \sum_{\mathbf{v} \in T^k} \phi\big( d(\mathbf{v}),
    X_\mathbf{v} \big) \Big],
\end{equation}
for a continuous bounded $\phi \colon \mathscr{D}_k \times E^k \to \R_+$,
where $\mathscr{D}_k$ is the set of $(k+1)\times (k+1)$ distance
matrices, and where for $\mathbf{v} = (v_1, \dots, v_k)$ we have defined
\[
    X_\mathbf{v} = (X_{v_i})_{1 \le i \le k},\qquad 
    d(\mathbf{v}) = (d(v_i, v_j))_{0 \le i,j \le k},\quad v_0 = \rho.
\]
Roughly speaking, the moments encapsulate the law of the (rooted) subtree spanned
by $k$ typical vertices in the tree. The method of moments states that,
under a mild technical condition, convergence of the moments for $k \ge
1$ is sufficient to deduce vague convergence of the law of the tree in
the Gromov-weak topology. These notions are made precise in
Section~\ref{S:mmmSpaces}, where the definitions of the appropriate
topologies are recalled. This section also contains some adaptations of
known results on random metric spaces to the specific features of our
problem  (infinite measures and vague convergence, marks, Gromov-vague
topology).

The second ingredient of our proofs is a new many-to-few formula to
compute these moments. 
It provides an expression for \eqref{eq:momentIntro} in terms of a
Markov process indexed by the vertices of a ``uniform tree'' with $k$
leaves. In our application, in between branch points, this Markov process
reduces to a Markov chain $(\zeta_n)_{n \ge 0}$ defined as the
$h$-transform of the mean semi-group of the branching process, namely
\[
    \forall x \in E,\quad 
    \E\big[ f(\zeta_1) \mid \zeta_0 = x\big]
    =
    \frac{\E\big[ \angle{\Xi_x, hf} \big]}{h(x)},
\]
for any measurable $f \colon E \to \R_+$. (This process arises in
several formulations of spinal decomposition theorems and is sometimes
referred to as the \emph{spine process}.) Upon reaching a branch point of
degree $d$, the chain splits in $d$ daughters whose random locations
depend on the $d$-th moment of the point process $\Xi_x$, described in 
Section~\ref{SS:constructionTree}. Let $\tau$ be a tree and 
$(X^\tau_v)_{v \in \tau}$ be this Markov process indexed by $\tau$. Let
$\tau_\mathbf{v}$ be the subtree spanned by $\mathbf{v} \in T^k$. The
many-to-few formula establishes that 
\[
    \E_x\bigg[ \sum_{\substack{\mathbf{v} \in T^k\\ v_1 < \dots < v_k}} 
        \phi\big( \tau_{\mathbf{v}}, (X_w)_{w \in \tau_{\mathbf{v}}} \big) 
    \bigg]
    =
    \sum_{\tau \in \mathscr{T}_k} \E_x\Big[ 
        (\Delta_k \phi)\big( \tau, (X^\tau_w)_{w \in \tau} \big)
    \Big],
\]
where $\phi$ is a test function and $\Delta_k$ is a bias term defined in
\eqref{eq:delta} that only depends on the spatial locations at the leaves
and branch points of $\tau$. For the previous expression to make sense,
$T$ will be endowed with a natural planar ordering of its vertices and
the sum on the right-hand side runs over the set $\mathscr{T}_k$ of
planar rooted trees with $k$ labeled leaves. The variables $(X^\tau_v)_{v
\in \tau}$ and $\Delta_k$ are constructed formally in Section~\ref{S:moments}, 
where the many-to-few formula (Theorem~\ref{thm:manyToFew}) is proved.

Finally, in Section~\ref{S:proof}, we use Assumption~\ref{ass:criticality} 
and the many-to-few formula to compute the limit of the moments of our
class of spatial branching processes. By the method of moments, this
will prove Theorem~\ref{thm:GromovVague}. The rest of the section is dedicated to deriving
tightness in the Gromov--Hausdorff--Prohorov topology from
Assumption~\ref{ass:kolmogorov}, which will prove Theorem~\ref{thm:GHP}. The key is
to use an adaptation of the tightness result derived in \cite{athreya16}.

\section{Preliminaries on planar trees}

\subsection{Discrete planar trees}

A rooted planar tree is a tree with a distinguished vertex and a total
order on the offspring of each vertex. It is standard to encode rooted
planar trees as subsets of the set of finite words
\[
    \mathcal{V} \coloneqq \{ \varnothing \} \cup \bigcup_{n=1}^\infty
    \N^n.
\]
As usual, we write $vw$ for the concatenation of $v$ and $w$, and
$\abs{v}$ for the length (or the \emph{generation}) of $v$. The set
$\mathcal{V}$ is naturally endowed with the lexicographical order $\le$,
and with a partial order $\preceq$ such that $v \preceq w$ if $v$
is an ancestor of $w$, that is if $w = vv'$ for some $v' \in
\mathcal{V}$. We write $v \wedge w$ for the most-recent common ancestor
of $v$ and $w$, that is the longest word $u$ such that $u \preceq v$ and
$u \preceq w$. 

A subset $\tau \subseteq \mathcal{V}$ is a \emph{planar tree} if
\begin{enumerate}
    \item[(i)] for every $v \in \tau$, if $w \preceq v$ then $w \in
        \tau$; and
    \item[(ii)] for every $v \in \mathcal{V}$ there is an integer $d_v(\tau)$
        such that 
        \[
            vi \in \tau \iff i \le d_v(\tau).
        \]
\end{enumerate}
The integer $d_v(\tau)$ is called the \emph{(out-)degree} of $v$ in
$\tau$. We introduce 
\[
    \mathcal{L}(\tau) = \{v \in \mathcal{V} : d_v(\tau) = 0 \},
    \qquad
    \mathcal{B}(\tau) = \{v \in \mathcal{V} : d_v(\tau) > 1 \},
\]
the set of \emph{leaves} and \emph{branch points} of
$\tau$. We denote by $\mathscr{T}_k$ the set of planar trees with $k$
leaves

Fix a vector of $k$ vertices $\mathbf{v} = (v_1, \dots, v_k) \in
\mathcal{V}^k$. We will denote by 
\[
    \tau_{\mathbf{v}} = \bigcup_{i=1}^k \{ w \in \mathcal{V} : w \preceq v_i \}
\]
the subtree spanned by $\mathbf{v}$. Note that, strictly speaking,
$\tau_{\mathbf{v}}$ is not a planar tree as defined above since the
children of a given vertex might not be indexed by consecutive integers
and hence might not fulfill (ii). However, there is a unique way to
relabel the elements of $\tau_{\mathbf{v}}$ to a planar tree that
preserves for each vertex its degree and the ordering of its children. In
practice, we do not make a distinction between $\tau_{\mathbf{v}}$ and
this relabeling.

\subsection{Two constructions of discrete planar trees}
\label{SS:constructionTree}

In this section, we provide two ways to generate a planar trees with $k$
leaves. The first one is a recursive construction obtained by decomposing
a tree into the branch leading to the first branch point and the subtrees
attached to it. It will allow us to obtain recursive formulas for the
moments of branching processes in terms of lower order moments. The
second one is an encoding of a planar trees as elements of $\N^{2k-1}$,
obtained by recording the heights of successive leaves and branch points.
It is this second construction that will allow us to take a limit to
continuous trees.

\begin{figure}[t]
    \centering
    \includegraphics[width=\textwidth]{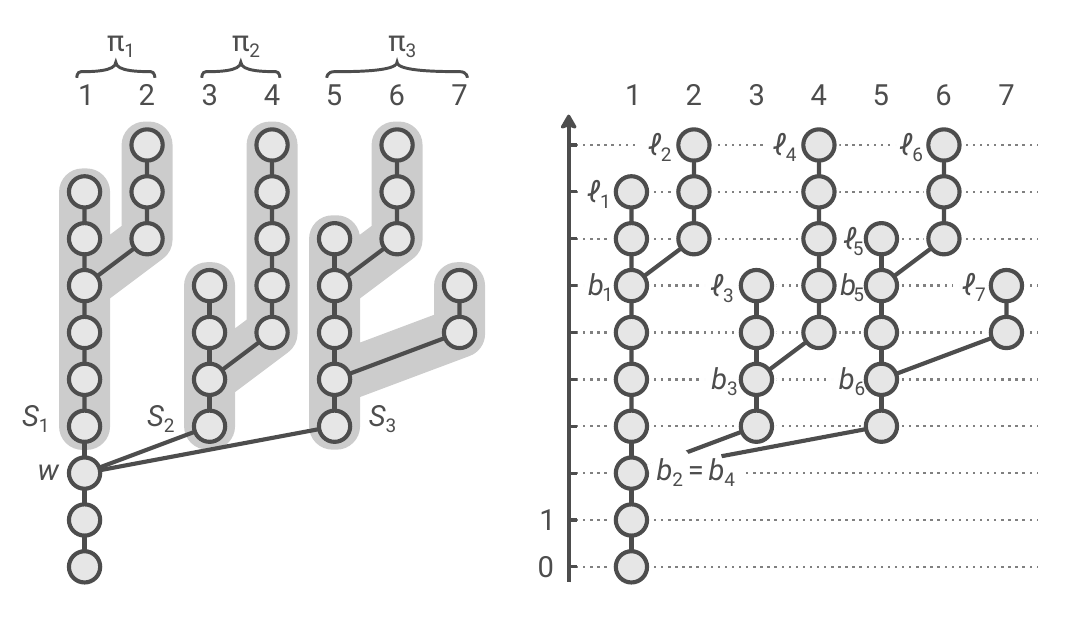}
    \caption{\emph{Left}: A tree is decomposed at its first branch point
    into three subtrees $S_1$, $S_2$, and $S_3$; the corresponding leaf
    partition is displayed on top. \emph{Right}: Point process
    construction of the same tree.}
    \label{fig:constructions}
\end{figure}

\paragraph{Recursive construction.} We will need a few definitions, which
are all illustrated in the left part of Figure~\ref{fig:constructions}.
Fix some tree $\tau \in \mathscr{T}_k$ with leaves $v_1 < \dots < v_k$
and let $w = v_1 \wedge \dots \wedge v_k$ be the first branch point. We
construct a partition $\pi(\tau)$ as 
\[
    i \sim_{\pi(\tau)} j \iff v_i \wedge v_j \ne w.
\]
In words, removing $w$ splits $\tau$ into smaller children subtrees, and
$i \sim_{\pi(\tau)} j$ if and only if the leaves $v_i$ and $v_j$ belong
to the same subtree. We let $\abs{\pi(\tau)}$ denote the number of blocks
of the partition, and $(\pi_i(\tau))_i$ be these blocks, ordered
according to their least element. Note that, due to the order structure,
the blocks of $\pi(\tau)$ are made of consecutive integers. Equivalently,
we could have recorded the vector of lengths of those blocks as $c =
(\abs{\pi_i(\tau)})_i$, which is a vector summing to $k$ referred to as a
\emph{composition of $k$}. Finally, for $i \le d_w(\tau)$, let us define
$S_i(\tau)$ as the $i$-th subtree attached to $w$. Namely,
\[
    \forall i \le d_w(\tau),\quad S_i(\tau) = \{ v \in \mathcal{V} : (w,i,v) \in \tau \}.
\]
It should hopefully be clear that the map
\[
    \tau \mapsto \big(\abs{w}, \big( S_i(\tau) \big)_{i \le \abs{\pi(\tau)}} \big)
\]
is a bijection from $\mathscr{T}_k$ to the set
\[
    \bigcup_{\pi} \N \times 
    \Big( \mathscr{T}_{\abs{\pi_1}} 
    \times \dots \times 
    \mathscr{T}_{\abs{\pi_{\abs{\pi}}}} \Big),
\]
where the union runs over all partitions of $[k]$ with blocks made of
consecutive integers, except the partition into singletons $\{ \{i\}, i
\le k\}$. A short proof is given in Lemma~\ref{lem:inductiveConstruction}.

\paragraph{Point process construction.} Fix a tree $\tau \in
\mathscr{T}_k$ with leaves $v_1 < \dots < v_k$. We define 
\[
    \forall i \le k,\quad \ell_i(\tau) = \abs{v_i},\qquad 
    \forall i \le k-1,\quad b_i(\tau) = \abs{v_i \wedge v_{i+1}}.
\]
The vectors $\bm{\uell}(\tau) = (\ell_i(\tau))_{i \le k}$ and
$\mathbf{b}(\tau) = (b_i(\tau))_{i \le k-1}$ encode, respectively, the
heights of the leaves and the heights of the branch points between
consecutive leaves. These are illustrated on the right part of
Figure~\ref{fig:constructions}. Again, it is not hard to see that the map
\[
    \tau \mapsto \big( \bm{\uell}(\tau), \mathbf{b}(\tau) \big)
\]
is a bijection from $\mathscr{T}_k$ to
\[
    \Big\{ (\bm{\uell}, \mathbf{b}) \in \N^k \times \N^{k-1} : b_i <
    \ell_i \wedge \ell_{i+1}, \, i \le k-1 \Big\}.
\]
We provide a short proof of this fact in Lemma~\ref{lem:cppConstruction}.
The inverse application is obtained starting from a branch of length
$\ell_1 +1$, and iteratively gluing a branch of length $\ell_{i+1}-b_i$
at the right-most vertex of the tree with height $b_i$, see Figure~\ref{fig:constructions}.

An interesting consequence of this construction is that the set of all
planar trees with $k$ leaves at height $\ell$, called ultrametric
trees, has a natural bijection with $\{0, \dots, \ell-1\}^{k-1}$. This
observation is at the heart of the \emph{coalescent point process}
construction of the genealogy of branching processes \cite{Popovic2004,
Lambert2013}, and more generally of the \emph{comb} encoding of
ultrametric spaces \cite{Lambert2017, foutelrodier2021}.

\subsection{Continuous trees, distance matrices}

\paragraph{Continuous trees.} The point process 
encoding of planar trees allows us to embed $\mathscr{T}_k$ into
a subset of $\N^{2k-1}$. By extending this construction to $\R_+^{2k-1}$,
we obtain a convenient notion of continuous trees with $k$ leaves and 
a topology for these objects. More precisely, define a continuous planar
tree as an element $\theta \in \mathscr{T}_{k,\mathrm{c}}$, where
\[
    \mathscr{T}_{k,\mathrm{c}} \coloneqq 
    \big\{  
        (\bm{\uell}, \mathbf{b}) \in \R_+^k \times \R_+^{k-1} : b_i <
        \ell_i \wedge \ell_{i+1},\, i \le k-1 
    \big\}.
\]
We write $\bm{\uell}(\theta)$ and $\mathbf{b}(\theta)$ for the
leaves and branch points heights of $\theta = (\bm{\uell}, \mathbf{b})$.
As in the discrete case, the tree $\theta = (\bm{\uell}, \mathbf{b})$ can
be obtained by starting from an initial branch corresponding to the
interval $[0, \ell_1]$ and inductively gluing an interval $[b_i,
\ell_{i+1}]$ on the right-most vertex at height $b_i$ of the tree
constructed at step $i$. We say that a sequence of planar trees
$(\theta_n)_{n \ge 1}$ converges to $\theta$ if 
\[
    \big( \bm{\uell}(\theta_n), \mathbf{b}(\theta_n) \big) \goesto{n \to \infty} 
    \big( \bm{\uell}(\theta), \mathbf{b}(\theta) \big).
\]
Moreover, if $a > 0$, we let $a \theta = (a\bm{\uell}(\theta),
a\mathbf{b}(\theta))$ be the tree obtained by rescaling the branch lengths of
$\theta$ by $a$. Finally, $\mathscr{T}_{k,\mathrm{c}}$ is endowed with 
two natural measures. The first one is a ``uniform'' measure $\Lambda_k$ defined as
\begin{equation} \label{eq:uniformTree}
    \int_{\mathscr{T}_{k,\mathrm{c}}} f( \theta ) \Lambda_k(\diff \theta)
    =
    \int_{\R_+^k} \diff \bm{\uell} \int_{\R_+^{k-1}} \diff \mathbf{b}
    \prod_{i=1}^{k-1} \indic_{\{ b_i < \ell_i \wedge \ell_{i+1} \}}
    f( \bm{\uell}, \mathbf{b} ).
\end{equation}
Note that this measure is infinite and only puts mass on binary trees.
The second one is a uniform measure $\widetilde{\Lambda}_k$ on the set of
ultrametric trees with $k$ leaves at height $1$, defined as 
\begin{equation} \label{eq:uniformUltrametric}
    \int_{\mathscr{T}_{k,\mathrm{c}}} f( \theta ) \widetilde{\Lambda}_k(\diff \theta)
    =
    \int_{[0,1]^{k-1}} f( \mathbf{1}_k, \mathbf{b} ) \,\diff \mathbf{b},
\end{equation}
where $\mathbf{1}_k = (1,\dots,1)$.

\paragraph{Distance matrices.} The (pointed) Gromov-weak topology relies
on studying the metric structure of a sample of size $k$ chosen uniformly
from a metric measure space. In its usual formulation, this metric
structure is given by a $(k+1) \times (k+1)$ matrix of pairwise
distances between the vector made of the $k$ points and the root. Let us
denote by $\mathscr{D}_k$ the set of $(k+1) \times (k+1)$ distance
matrices:
\begin{equation} \label{eq:defSetDistanceMatrix}
    \mathscr{D}_k = \{ (d_{ij})_{0 \le i,j \le k} : d_{ii} = 0,\, d_{ij} =
    d_{ji},\, d_{ij} \le d_{ip} + d_{pj},\, 1 \le i,j,p \le k \}.
\end{equation}
Since all the metric spaces that we consider in this work are trees, it
more natural to encode a sample of size $k$ as an element of
$\mathscr{T}_{k,\mathrm{c}}$ rather than as a distance matrix. The two
points of view are of course equivalent. To make the connexion explicit,
we let $D(\theta)$ be the $(k+1)\times(k+1)$ matrix of pairwise distances
between the leaves and the root of $\theta \in \mathscr{T}_{k,\mathrm{c}}$,
\begin{equation} \label{eq:defDistanceMatrix}
    \forall 0 \le i < j \le k,\quad
    D_{ji}(\theta) = D_{ij}(\theta)
    =
    \begin{dcases}
        \ell_j &\text{if $i = 0$,} \\
        \ell_i + \ell_j - \min \{ b_i, \dots, b_{j-1} \} 
               &\text{if $i \ge 1$},
    \end{dcases}
\end{equation}
and $D_{ii}(\theta) = 0$.

\subsection{Marked trees}

The class of branching processes that we are
interested in defines a planar tree structure augmented with types (that
we also call marks or locations) on the vertices. Formally, if $E$ is a
fixed mark space we define a (discrete) \emph{marked tree} as a pair
$\tau^* = (\tau, (x_v)_{v \in \tau})$, where $\tau$ is a planar tree 
and $x_v \in E$ for each $v \in \tau$. We let $\mathscr{T}_k^*$ be the
set of all discrete marked trees with $k$ leaves. Without further
mention, we extend all of the operations that we have defined on trees to
their obvious counterparts for marked trees.

In the continuous setting, the natural extension of this definition would
be that a marked tree is a pair $(\theta, (x_v)_{v \in \theta})$ where
the second coordinate is a collection of marks indexed by $\theta$.
Although this notion is intuitive, giving a formal definition would
rapidly become cumbersome. For our purpose we will only need to record
the types of the $k$ leaves of $\theta$. Accordingly, we define a
continuous marked tree as a pair $\theta^* = (\theta, (x_i)_{i \le k})$,
with the interpretation that $x_i$ is the type of the $i$-th leave of
$\theta$, and let $\mathscr{T}^*_{\mathrm{c},k} = \mathscr{T}_k \times
E^k$. We also let $\mathscr{D}_k^* = \mathscr{D}_k \times E^k$ 
be the set of marked distance matrices.

\section{The moments of branching processes}
\label{S:moments}

\subsection{Notation and tree-index Markov chain}

\paragraph{Construction of the branching process.}
We give a brief construction of the spatial branching processes that are
described in the introduction. Recall that $\Xi_x$ stands for the point
process encoding the offspring of a particle at $x \in E$. We
construct iteratively a tree $T_n$ corresponding to the first $n$
generations of the branching process as follows:
\begin{itemize}
    \item Let $T_0 = \{ \varnothing \}$ and $X_\varnothing = x_0$.
    \item For $n \ge 0$, conditional on $(T_n, (X_v)_{v \in T_n})$,
        let $(\Xi_v)_{v \in T_n, \abs{v} = n}$ be independent point
        processes such that $\Xi_v = \sum_i \xi_{v,i} \sim \Xi_{X_v}$.
        Define
        \[
            T_{n+1} = T_n \cup \bigcup_{v \in T_n, \abs{v} = n} 
            \big\{ vi : i \le \abs{\Xi_v} \big\},
            \qquad
            X_{vi} = \xi_{v, i}.
        \]
    \item Let $T = \cup_{n \ge 0} T_n$ and $T^* = (T, (X_v)_{v \in T})$.
\end{itemize}

Recall that $\P_x$ is such that $T^*$ starts from a single ancestor at
$X_\varnothing = x$ under $\P_x$. It will considerably ease the notation
to let $\Xi$ be a generic point process such that, under $\P_x$, $\Xi$ is
distributed as $\Xi_x$. We will write its atoms as
\[
    \Xi = \sum_{i=1}^{\abs{\Xi}} \delta_{\xi_i}.
\]

\paragraph{Biasing function $\psi$ and spine.}
All our results in this section depend on fixing a positive function
$\psi \colon E \to (0, \infty)$ that ``weights'' the contribution of
individuals in the population according to their types, in the spirit
of \cite{Bansaye2024, Marguet2019}. As we shall see
later, two important examples that will lead to simplifications are $\psi
\equiv 1$, which provides a Feynman--Kac representation of the moments,
and $\psi \equiv h$ where $h$ is a harmonic function of the dynamics.

Suppose that $\psi$ has been fixed. Define the $d$-th factorial
$\psi$-moment of the offspring law as 
\begin{equation} \label{eq:defMoment}
    \forall d \ge 1,\: \forall x \in E,\quad m^\psi_d(x) = 
    \E_x\Big[ \sum_{\substack{i_1, \dots, i_d = 1\\\text{$(i_p)_p$ distinct}}}^{\abs{\Xi}}
    \prod_{i=1}^d \psi(\xi_i) \Big],
\end{equation}
and, provided that $m_d^\psi(x) < \infty$, let $(\chi_{i,d})_{i \le d}$
be the random variables with distribution
\begin{equation} \label{eq:jumpDistr}
    \forall x \in E,\quad \E_x\big[ f( \chi_{1,d}, \dots, \chi_{d,d}) \big] = \frac{1}{m^\psi_d(x)} 
    \E_x\Big[ \sum_{\substack{i_1, \dots, i_d = 1\\\text{$(i_p)_p$ distinct}}}^{\abs{\Xi}}
    \prod_{i=1}^d \psi(\xi_i) f(\xi_1,\dots, \xi_d) \Big].
\end{equation}
In words, these random variables are obtained by first biasing the
offspring distribution by its $d$-th $\psi$-moment and then sampling $d$
particles, each with probability proportional to $\psi$.

The case $d = 1$ will play an essential role. We let $(\zeta_n)_{n \ge
0}$ denote the Markov chain with transition probability
\[
    \E[ f(\zeta_1) \mid \zeta_0 = x ] = \E_x[ f(\chi_{1,1}) ]
    = \frac{1}{m^\psi_1(x)} 
    \E_x\Big[ \sum_{i=1}^{\abs{\Xi}} \psi(\xi_i) f(\xi_i) \Big].
\]
We call $(\zeta_n)_{n \ge 1}$ the \emph{spinal Markov chain} (or spinal
process). It is the usual Markov chain that arises in the study of many
classes of branching process, through a many-to-one formula or a spinal
decomposition result \cite{kurtz1997conceptual, biggins2004measure,
shi2015branching}. Let us introduce a final notation
\begin{equation*} 
    \forall x \in E,\quad \lambda^\psi(x) = \frac{\E_x[\angle{\Xi,
    \psi}]}{\psi(x)} = \frac{m_1^\psi(x)}{\psi(x)},
\end{equation*}
which gives the expected growth of the population size, starting from one
particle at $x$ and when size is measured according to $\psi$.

\paragraph{Tree-indexed Markov chain.} Our final task before stating the
many-to-few formula is to define the notion of a Markov chain indexed by
a tree. It is a simple procedure to assign a collection of types
$(X^\tau_v)_{v \in \tau}$ to the vertices of a fixed tree $\tau \in
\mathscr{T}_k$, which is described for instance in \cite{Marckert2008}.
For $\tau \in \mathscr{T}_k$, define $(X^\tau_v)_{v \in \tau}$ as
follows:
\begin{itemize}
    \item Start from $X^\tau_\varnothing = x$.
    \item Conditional on $(X^\tau_v)_{\abs{v} \le n}$, the variables
        $(X^\tau_{vi})_{i \le d_v}$ are distributed as
        \[
            (X^\tau_{v1}, \dots, X^\tau_{vd_v}) \overset{\mathrm{(d)}}{=} 
            (\chi_{1,d_v}, \dots, \chi_{d_v, d_v}) \: \text{under $\P_{X_v}$},
        \]
        independently for different vertices $v$.
\end{itemize}
We denote by $\Q^\psi_{x,\tau}$ the distribution of the random marked tree
$(\tau, (X^\tau_v)_{v \in \tau})$, which is a probability measure on
$\mathscr{T}^*_k$.

The trees that we consider will typically have a large height of order
$n$ but a fixed number of leaves $k$, and thus most of their vertices
have degree $1$. The image that one should keep in mind is that the
process indexed by such a tree is obtained by running independent copies
of the spinal process $(\zeta_n)_{n \ge 0}$ along the branches of the
tree, allocating the offspring types as in \eqref{eq:jumpDistr} whenever
a branch point is reached.

\subsection{The many-to-few formula}

Fix $k \ge 1$. We define the $k$-th \emph{planar (factorial) moment} of the
branching process as the measure $\M^k_x$ on the space $\mathscr{T}_k^*$
of marked trees with $k$ leaves such that 
\begin{equation} \label{eq:momentMeasure}
    \M^{k}_x[F] = \E_x\bigg[
        \sum_{\substack{\mathbf{v} \in T^k\\v_1 < \dots < v_k}}
        F \big(\tau^*_\mathbf{v} \big) 
    \bigg],
\end{equation}
for any non-negative map $F \colon \mathscr{T}^*_k \to \R_+$ and where
$\tau^*_\mathbf{v} = (\tau_\mathbf{v}, (X_w)_{w \in \tau_\mathbf{v}})$ is
the marked subtree spanned by $\tau_{\mathbf{v}}$.
Note that $\tau_{\mathbf{v}}$ might fail to have $k$ leaves. In that case
we use the convention that $F(\tau^*_{\mathbf{v}}) = 0$, which amounts to
removing the corresponding terms from the sum.

We can now state the main result of this section, which gives an
expression for $\M^k_x$ in terms of the Markov chain indexed by a
``uniform'' element of $\mathscr{T}_k$ and an appropriate biasing term
defined as
\begin{equation} \label{eq:delta}
    \forall \tau^* \in \mathscr{T}^*_k,\quad
    \Delta^\psi_k(\tau^*) = 
    \prod_{v \in \tau} \lambda^\psi(x_v)
    \prod_{v \in \mathcal{B}(\tau)} \frac{m^\psi_{d_v}(x_v)}{d_v!\,\psi(x_v)\lambda^\psi(x_v)} 
    \prod_{v \in \mathcal{L}(\tau)} \frac{1}{\psi(x_v)\lambda^\psi(x_v)}.
\end{equation}
Recall that $\Q^\psi_{x,\tau}$ denotes the distribution of the pair
$(\tau, (X^\tau_v)_{v \in \tau})$ started from $x$, where $(X^\tau_v)_{v
\in \tau}$ is Markov process indexed by $\tau$.

\begin{theorem}[Many-to-few at all times] \label{thm:manyToFew}
    Fix some functional $\psi \colon E \to (0, \infty)$. Then, 
    \[
        \M^{k}_x[F] = \psi(x) \sum_{\tau \in \mathscr{T}_k} 
        \Q^\psi_{x,\tau}[\Delta^\psi_k F],
    \]
    for any non-negative functional $F \colon \mathscr{T}^*_k \to \R_+$.
\end{theorem}

\begin{remark}
    The idea of using trees to compute the moments of branching processes
    appears as early as in the work of Dynkin \cite{Dynkin1991} in the
    superprocess literature, see also \cite[Section~2.1]{etheridge2000introduction}.
    Our many-to-few formula is very similar in spirit. An important
    distinction, however, is that the trees appearing here provide
    information about the genealogy of the branching process, whereas they
    are only used as a formal tool to compute the moments in \cite{Dynkin1991,
    etheridge2000introduction}.
\end{remark}

Another interpretation of this result is that the bias term
$\Delta^\psi_k$ is the Radon--Nikodym derivative of $\M^k_x$ with respect
to the infinite measure $\sum_{\tau \in \mathscr{T}_k} \Q^\psi_{x,\tau}$.
The latter measure is the Markov chain indexed by a ``uniform tree'' with
$k$ leaves.

This result is proved in the following section. As an intermediate step,
we will obtain an equivalent formulation of the many-to-few formula
(Proposition~\ref{prop:many2fewRecursion}) which can be of independent
interest. It provides a recursive expression of the moment measure
$\M^k_x$ in terms of lower order moment measures. Before moving to the
proofs, let us give a few important examples of choices for $\psi$.

\paragraph{Harmonic function.} Suppose that there exists $\lambda > 0$
(fixed) and $h \colon E \to (0, \infty)$ such that 
\[
    \forall x \in E,\quad \E_x[ \angle{\Xi, h} ] = \lambda h(x).
\]
The pair $(\lambda, h)$ can be seen as a solution to an eigenproblem for
the mean semi-group of the branching process and $h$ is a \emph{harmonic
function}. In that case, letting $\psi \equiv h$, we see that
$\lambda^\psi(x) \equiv \lambda$ and the biasing term \eqref{eq:delta}
simply becomes
\[
    \forall \tau^* \in \mathscr{T}^*_k,\quad 
    \Delta_k(\tau^*) = 
    \lambda^{\abs{\tau}-k-b}
    \prod_{v \in \mathcal{B}(\tau)} \frac{m_{d_v}(x_v)}{d_v!\,h(x_v)} 
    \prod_{v \in \mathcal{L}(\tau)} \frac{1}{h(x_v)}.
\]
The corresponding spine Markov chain $(\zeta_n)_{n \ge 0}$ is Doob's
harmonic transform of the mean semi-group of the branching process.

\paragraph{Feynman--Kac formula.} If $\psi \equiv 1$, $(\zeta_n)_{n \ge
0}$ is simply obtained by size-biasing the offspring point process and
jumping to the location of a uniformly chosen particle. The many-to-one
formula (for $k=1$) then becomes
\[
    \E_x\Big[ \sum_{\abs{v} = n} f(X_v) \Big]
    = 
    \E_x\Big[ f(\zeta_n) \prod_{m=0}^{n-1} \E_{\zeta_m}[\abs{\Xi}] \Big],
\]
which is the usual Feynman--Kac formula for the mean semi-group of the
branching process. For $k \ge 1$, the bias term becomes
\[
    \Delta_k(\tau^*) = 
    \prod_{\substack{v \in \tau\\ v \notin \mathcal{L}(\tau)}} 
    \frac{1}{d_v!} \E_{x_v}\big[ \abs{\Xi}^{(d_v)} \big],
\]
where $n^{(d)} = n(n-1)\dots (n-k+1)$. It only depends on the
factorial moments of the offspring point process.

\subsection{Proof of the many-to-few formula}

It will be more convenient to work with the factorial moment
of the process biased by $\psi$, which we define as the measure
$\M^{k,\psi}_x$ on $\mathscr{T}^*_k$ such that, for any $F \colon
\mathscr{T}^*_k \to \R_+$,
\[
    \M^{k,\psi}_x[F] = \E_x\Big[
        \sum_{\substack{\mathbf{v} \in T^k\\v_1 < \dots < v_k}}
        F \big(\tau^*_\mathbf{v} \big) 
        \prod_{i=1}^k \psi( X_{v_i} )
    \Big].
\]

\subsubsection{Many-to-one formula}

We start by deriving an expression for the case $k=1$. In that case,
$\M^{1,\psi}_x$ can be identified as a measure on the space of (discrete) paths
$\cup_{n \ge 0} E^{n+1}$.

\begin{proposition}[Many-to-one at all times] \label{prop:manyToOne}
    For any functional $F \colon \cup_{n \ge 0} E^{n+1} \to \R_+$,
    \[
        \forall x \in E,\quad \M^{1,\psi}_x[F] = 
        \psi(x) \sum_{n \ge 0} \E_x\Big[ 
            F\big( (\zeta_i)_{i\le n} \big)
            \prod_{i=0}^{n-1} \lambda^\psi(\zeta_i) 
        \Big].
    \]
\end{proposition}

\begin{proof}
    It is sufficient to show that the identity holds for functionals $F$
    whose support in included in the set of paths with length smaller
    than $n$, for any $n \ge 0$. We prove the result by induction on this
    length $n \ge 0$. For $n = 0$, by definition
    \[
        \M^{1,\psi}_x[F] 
        = \E_x\big[ F( X_\varnothing ) \psi(X_\varnothing) \big]
        = \psi(x) F(x).
    \]
    For $n \ge 1$, applying the branching property after the first
    reproduction event,
    \begin{align*}
        \M^{1,\psi}_x[F] 
        &= \E_x\Big[ 
            \sum_{v \in T} F\big((X_w)_{w \preceq v}\big) \psi(X_v)
        \Big] \\
        &= \psi(x)F(x) +
        \E_x\Big[ \sum_{i=1}^{\abs{\Xi}}
            \M^{1,\psi}_{\xi_{i}}\big[F\big( (x, \cdot) \big)\big]
        \Big] \\
        &\overset{(\dagger)}{=} \psi(x)F(x) +
        \E_x\Big[ \sum_{i=1}^{\abs{\Xi}}
            \sum_{n \ge 0} \psi(\xi_i)
            \E_{\xi_i}\Big[ F\big( (x, \zeta_j)_{j \le n}\big)
            \prod_{j=0}^{n-1} \lambda^\psi(\zeta_j) \Big]
        \Big] \\
        &\overset{(\ddagger)}{=} \psi(x)F(x) + 
        \sum_{n \ge 0} 
        \E_x\big[ \angle{\Xi, \psi} \big]
        \E_x\Big[ 
            \E_{\zeta_1} \Big[ F\big( (x, \zeta_j)_{j \le n}\big)
            \prod_{j=0}^{n-1} \lambda^\psi(\zeta_j) \Big]
        \Big] \\
        &= \psi(x)
        \bigg( F(x) + \lambda^\psi(x) 
        \sum_{n \ge 0} 
        \E_x\Big[ 
            F\big((\zeta_j)_{j \le n+1}\big)
            \prod_{j=1}^n \lambda^\psi(\zeta_j)
        \Big] \bigg) \\
        &=
        \psi(x) \sum_{n \ge 0} \E_x\Big[ 
            F\big( (\zeta_j)_{j\le n} \big)
            \prod_{j=0}^{n-1} \lambda^\psi(\zeta_j) 
        \Big],
    \end{align*}
    where we have used our induction in $(\dagger)$ and
    the definition of $(\zeta_n)_{n \ge 0}$ in $(\ddagger)$.
\end{proof}

\begin{remark}
    One could obtain a more precise result. Note that
    \[
        \forall n \ge 0,\quad W^\psi_n = \sum_{v \in T, \abs{v} = n}
        \psi(X_v) \prod_{w \prec v} \lambda^\psi(X_w)
    \]
    is a martingale. Standard arguments would show that the martingale
    change of measure $W_n \diff \P_x$ can be represented as the law of a
    spine with distribution $(\zeta_n)_{n \ge 0}$ on which are grafted
    independent subtrees distributed as the original process. The
    many-to-few formula is simply the projection of this result on the law
    of the spine, see \cite[Chapter~4]{shi2015branching} for instance.
\end{remark}

\subsubsection{Many-to-few formula: factorized form}

We establish a first version of our many-to-few formula which provides an
inductive expression for the $k$-th moment measure in terms of lower
order moments. It relies on decomposing a marked tree $\tau^* \in
\mathscr{T}^*_k$ at its first branch point. Recall the notation of
Section~\ref{SS:constructionTree}, in particular that $(S_i(\tau^*))_i$
are the marked subtrees obtained by removing the first branch point of
$\tau$, and that $\pi(\tau)$ is the corresponding partition of the
leaves. We will evaluate the moment measure $\M^{k,\psi}_x$ at
\emph{product functionals} of the form
\begin{equation} \label{eq:functionalM2F}
    \forall \tau^* \in \mathscr{T}^*_k,\quad 
    F(\tau^*) = \indic_{\{ \pi(\tau) = \pi\}} F_0\big((x_v)_{v \preceq w} \big)
    \prod_{i=1}^{\abs{\pi}} F_i(\tau^*_i),
\end{equation}
for some partition $\pi$ of $[k]$, functionals $F_i \colon
\mathscr{T}^*_{\abs{\pi_i}} \to \R_+$, and $F_0 \colon \cup_{n \ge 0} E^{n+1}
\to \R_+$. (In this expression $w = v_1 \wedge \dots \wedge v_k$.)

\begin{proposition}[Many-to-few, factorized form] \label{prop:many2fewRecursion}
    Let $F$ be a functional of the form \eqref{eq:functionalM2F}. Then
    \begin{equation*}
        \M^{k,\psi}_x[F] = 
        \psi(x) \sum_{n = 0}^\infty 
        \E_x\bigg[ \Big( \prod_{m=0}^{n-1} \lambda^\psi(\zeta_m)\Big)
            \frac{F_0((\zeta_m)_{m \le n})}{\abs{\pi}!\, \psi(\zeta_n)}
            \E_{\zeta_n}\Big[ \sum_{\substack{r_1, \dots, r_{\abs{\pi}} = 1\\ \text{\normalfont$(r_i)$ distinct}}}^{\abs{\Xi}} 
                \prod_{i=1}^{\abs{\pi}} \M^{\abs{\pi_i},\psi}_{\xi_{r_i}}\big[ F_i \big]
            \Big] \bigg],
    \end{equation*}
    where $(\zeta_n)_{n \ge 0}$ is the spine Markov chain. 
\end{proposition}

\begin{proof}
    Fix some $v_1 < \dots < v_k$. If $\pi(\tau_{\mathbf{v}}) = \pi$ and
    $w = v_1 \wedge \dots \wedge v_k$, there exist $\abs{\pi}$ integers
    $r_1 < \dots < r_{\abs{\pi}}$ giving the labels of the children of
    $w$ that correspond to the roots of the subtrees $(S_i(\tau^*))_i$. More
    formally, $\pi(\tau_{\mathbf{v}}) = \pi$ if and only if there is (a
    unique) $w \in \mathcal{V}$ and $r_1 <
    \dots < r_{\abs{\pi}}$ such that 
    \[
        \forall i \le \abs{\pi},\: \forall j \in \pi_i,\quad
        (w, r_i) \preceq v_j.
    \]
    Thus, we can write 
    \[
        \indic_{\{\pi(\tau_\mathbf{v}) = \pi\}}
        = \sum_{w \in \mathcal{V}} \sum_{r_1 < \dots < r_{\abs{\pi}}} 
        \indic_{\{ \forall i \le \abs{\pi},\:
        \forall j \in \pi_i,\: (w, r_i) \preceq v_j \}}.
    \]
    Therefore,
    \begin{align*}
        \M^{k,\psi}_x[F]
        &= \E_x\Big[ \sum_{\substack{\mathbf{v} \in T\\v_1 < \dots <
        v_k}} F(\tau^*_\mathbf{v}) \prod_{j=1}^k \psi(X_{v_j}) \Big] \\
        &= \E_x\Big[ \sum_{w \in \mathcal{V}} F_0\big((X_v)_{v \preceq w} \big)
            \sum_{r_1 < \dots < r_{\abs{\pi}}}
            \sum_{\substack{\mathbf{v} \in T\\v_1 < \dots < v_k}}
            \prod_{i=1}^{\abs{\pi}} F_i(S_i(\tau^*_\mathbf{v})) \prod_{j
            \in \pi_i} \psi(X_{v_j}) \indic_{\{ (w, r_i) \preceq v_j \}} 
        \Big].
    \end{align*}
    Now, note that $(w,r_i) \preceq v_j$ if and only if there exists
    $v'_j \in \mathcal{V}$ such that $v_j = (w,r_i,v'_j)$, and
    that in this case $S_i(\tau^*_\mathbf{v})$ is the subtree spanned by 
    $\mathbf{v}'_i = (v'_j)_{j \in \pi_i}$. Therefore,
    \begin{equation*}
        \sum_{\substack{\mathbf{v} \in T\\v_1 < \dots < v_k}}
        \prod_{i=1}^{\abs{\pi}} F_i(S_i(\tau^*_{\mathbf{v}})) 
        \indic_{\{ \forall j \in \pi_i,\, (w, r_i) \preceq v_j \}}
        = 
        \sum_{\substack{\mathbf{v}' \in \mathcal{V}\\v_1 < \dots < v_k}} 
        \prod_{i=1}^{\abs{\pi}} F_i(\tau^*_{\mathbf{v}_i'}) 
        \indic_{\{ \forall j \in \pi_i,\, (w, r_i, v'_j) \in T\}}.
    \end{equation*}
    The branching property applied at vertices $(w,r_1), \dots,
    (w,r_{\abs{\pi}})$ yields that  
    \begin{align*}
        \M^{k,\psi}_x[F]
        &= \E_x\Big[ \sum_{w \in \mathcal{V}} F_0\big((X_v)_{v \preceq w} \big)
           \sum_{r_1 < \dots < r_{\abs{\pi}}}
           \prod_{i=1}^{\abs{\pi}} \indic_{\{ (w,r_i) \in T\}}
           \M^{\abs{\pi_i}, \psi}_{X_{(w,r_i)}}\big[ F_i \big]
           \Big] \\
        &= \E_x\Big[ \sum_{w \in \mathcal{V}} F_0\big((X_v)_{v \preceq w} \big)
            \indic_{\{ w \in T \}}
            \frac{1}{\abs{\pi}!} \E_{X_w}\Big[ 
                \sum_{\substack{r_1, \dots, r_{\abs{\pi}} = 1\\ \text{\normalfont$(r_i)$ distinct}}}^{\abs{\Xi}} 
                \prod_{i=1}^{\abs{\pi}} \M^{\abs{\pi_i},\psi}_{\xi_{r_i}}\big[ F_i \big]
            \Big] \Big].
    \end{align*}
    The proof is now ended by applying the many-to-one formula
    (Proposition~\ref{prop:manyToOne}) to the above expression, to obtain
    that
    \[
        \M^{k,\psi}_x[F]
        =
        \sum_{n \ge 0} \psi(x) \E_x\bigg[
            \Big( \prod_{m=0}^{n-1}\lambda^\psi(\zeta_m) \Big) \frac{F_0\big( (\zeta_m)_{m \le n} \big)}{\psi(\zeta_n) \abs{\pi}!}
            \E_{\zeta_n}\Big[ 
            \sum_{\substack{r_1, \dots, r_{\abs{\pi}} = 1\\ \text{\normalfont$(r_i)$ distinct}}}^{\abs{\Xi}} 
            \prod_{i=1}^{\abs{\pi}} \M^{\abs{\pi_i},\psi}_{\xi_{r_i}}\big[ F_i \big]
        \Big] \bigg]. \qedhere
    \]
\end{proof}

\subsubsection{Completing the proof}

\begin{proof}[Proof of Theorem~\ref{thm:manyToFew}]
    Fix a partition $\pi$, some trees $\tau_i \in \mathscr{T}_{\abs{\pi_i}}$ 
    for $i \le \abs{\pi}$. Consider a functional $F$ of the form
    \eqref{eq:functionalM2F}, where each $F_i$ is itself of the form
    \[
        \forall \tau^* \in \mathscr{T}^*_{\abs{\pi_i}}, \quad 
        F_i(\tau^*) = \tilde{F}_i(\tau^*) \indic_{\{ \tau = \tau_i\}}.
    \]
    Since $\M^{k,\psi}_x$ is obtained by biasing $\M^k_x$ by the product
    of the weights at the leaves, it is sufficient to prove that 
    \[
        \M^{k,\psi}_x[F] = \sum_{\tau \in \mathscr{T}_k}
        \Q^{\psi}_{x,\tau}[\widetilde{\Delta}_k \cdot F],
        \qquad
        \widetilde{\Delta}_k(\tau^*) = \Delta_k(\tau^*) 
        \prod_{v \in \mathcal{L}(\tau)} \psi(x_v).
    \]
    For $k = 1$, this identity is the content of the many-to-one formula,
    Proposition~\ref{prop:manyToOne}. By an induction on the number of
    leaves and Proposition~\ref{prop:many2fewRecursion},
    \begin{align}
        \M^{k,\psi}_x[F]
        &=
        \sum_{n \ge 0} \psi(x) \E_x\bigg[
            \Big( \prod_{m=0}^{n-1}\lambda^\psi(\zeta_m) \Big) 
            \frac{F_0\big( (\zeta_m)_{m \le n} \big)}{\psi(\zeta_n) \abs{\pi}!}
            \E_{\zeta_n}\Big[ 
            \sum_{\substack{r_1, \dots, r_{\abs{\pi}} = 1\\ \text{\normalfont$(r_i)$ distinct}}}^{\abs{\Xi}} 
            \prod_{i=1}^{\abs{\pi}} \M^{\abs{\pi_i},\psi}_{\xi_{r_i}}\big[ F_i \big]
        \Big] \bigg] \nonumber \\
        &\overset{(\dagger)}{=}
        \sum_{n \ge 0} \psi(x) \E_x\bigg[
            \Big( \prod_{m=0}^{n-1}\lambda^\psi(\zeta_m) \Big) 
            \frac{F_0\big( (\zeta_m)_{m \le n} \big) m^\psi_{\abs{\pi}}(\zeta_n)}{\psi(\zeta_n) 
            \abs{\pi}!} \E_{\zeta_n}\Big[ 
            \prod_{i=1}^{\abs{\pi}}
            \frac{\M^{\abs{\pi_i},\psi}_{\chi_{i,\abs{\pi}}}\big[ F_i
                \big]}{\psi(\chi_{i,\abs{\pi}})}
        \Big] \bigg] \nonumber \\
        &=
        \sum_{n \ge 0} \psi(x) \E_x\bigg[
            \Big( \prod_{m=0}^{n-1}\lambda^\psi(\zeta_m) \Big) 
            \frac{F_0\big( (\zeta_m)_{m \le n} \big) 
            m^\psi_{\abs{\pi}}(\zeta_n)}{\psi(\zeta_n) \abs{\pi}!}
            \E_{\zeta_n}\Big[ \prod_{i=1}^{\abs{\pi}} 
                \Q_{\chi_{i,\abs{\pi}}, \tau_i}^\psi \big[ \widetilde{\Delta}_{\abs{\pi_i}} F_i \big]
        \Big] \bigg], \label{eq:m2f1}
    \end{align}
    where ($\dagger$) follows by definition of $(\chi_{i,d})_{i \le d}$
    in \eqref{eq:jumpDistr}. The result now follows from two
    observations. First, if $G$ is a product functional, namely, 
    \begin{equation*} 
        \forall \tau^* \in \mathscr{T}^*_k,\quad 
        G(\tau^*) = \indic_{\{ \pi(\tau) = \pi\}} G_0\big((x_v)_{v \preceq w} \big)
        \prod_{i=1}^{\abs{\pi}} G_i(\tau^*_i),
    \end{equation*}
    by construction of the tree-indexed Markov chain and a Markov property
    at the first branch point of $\tau$ 
    \begin{equation} \label{eq:prodSnake}
        \Q^{\psi}_{x,\tau}\big[ G \big]
        =
        \E_x\bigg[ G_0\big( (\zeta_m)_{m \le \abs{w}} \big)
        \E_{\zeta_{\abs{w}}}\Big[ 
            \prod_{i=1}^{\abs{\pi}}
            \Q^{\psi}_{\chi_{i,\abs{\pi}}, \tau_i}[G_i]
        \Big]
        \bigg].
    \end{equation}
    Second, the bias $\widetilde{\Delta}_k$ is itself of the product form:
    \[
        \widetilde{\Delta}_k(\tau^*) = \prod_{v \prec w}  \lambda^\psi(x_v)
        \frac{m^\psi_{\abs{\pi}}(x_w)}{\abs{\pi}!\, \psi(x_w)}
        \prod_{i=1}^{\abs{\pi}} \widetilde{\Delta}_{\abs{\pi_i}}(\tau^*_i).
    \]
    Therefore, taking as product functional     
    \begin{align*}
        G(\tau^*)
        &= \Big( \prod_{v \prec w} \lambda^\psi(x_v) \Big)
        \frac{F_0\big( (x_v)_{v \preceq w} \big) m^\psi_{\abs{\pi}}(x_w)}{\abs{\pi}!\, \psi(x_w)}
        \prod_{i=1}^{\abs{\pi}}
        \widetilde{\Delta}_{\abs{\pi_i}}(\tau^*_i) F_i(\tau^*_i) \\
        &= \widetilde{\Delta}_k(\tau^*) 
        F_0\big( (x_v)_{v \preceq w} \big)  \prod_{i=1}^{\abs{\pi}}
        F_i(\tau^*_i) 
        = \widetilde{\Delta}_k(\tau^*) F(\tau^*)
    \end{align*}
    in \eqref{eq:prodSnake} and plugging the result into \eqref{eq:m2f1}
    yields that
    \[
        \M^{k,\psi}_x[F] = \psi(x) \Q_x^{\tau, \psi}\big[ \widetilde{\Delta}_k \cdot F \big].
    \]
    This concludes the proof.
\end{proof}

\section{Marked metric measure spaces}
\label{S:mmmSpaces}

\subsection{Gromov topologies}
\label{SS:GromovTopologies}

Following \cite{greven2009convergence, Depperschmidt2011, athreya16},
a pointed marked metric measure space (mmm-space) is a quadruple
$\mathcal{X} = (X,d,\rho,\mu)$ where $(X,d)$ is a complete separable
metric space, $\rho \in X$ is a pointed element (the root) and $\mu$ is a
finite measure on $X \times E$. We let $\mathscr{X}$ denote
the set of equivalence classes of mmm-spaces, where $\mathcal{X}$ and
$\mathcal{X}'$ are equivalent if there exists a bijective isometry
\[
    \iota \colon \{ \rho \} \cup \supp_X \mu \to \{ \rho' \} \cup
    \supp_{X'} \mu'
\]
such that $\iota(\rho) = \rho'$ and $\mu \circ \tilde{\iota}^{-1} =
\mu'$, where $\tilde{\iota}(x,e) = (\iota(x), e)$. Here we have used the
notation $\supp_X \mu$ for the support of the projection of $\mu$ on its
first coordinate $X$. For a Borel set $A \subseteq X$ we will write
$\abs{A} = \mu(A \times E)$ for its mass. We briefly recall the
definition of three topologies on $\mathscr{X}$ that we will need. Recall
that $\mathscr{D}^*_k$ stands for the set of marked distance matrices
defined in \eqref{eq:defSetDistanceMatrix}. 

\paragraph{Gromov-weak topology.}
For $k \ge 1$ and $\phi \colon \mathscr{D}^*_k \to \R$, define
a functional $\Phi \colon \mathscr{X} \to \R$ as
\begin{equation} \label{eq:polynomial}
    \Phi(\mathcal{X}) = 
    \int_{(X\times E)^k} \phi\big( d(\mathbf{u}), \mathbf{e} \big) 
    \mu^{\otimes k}(\diff \mathbf{u}, \diff \mathbf{e}),
\end{equation}
where
\[
    \forall \mathbf{u} = (u_i)_{1 \le i \le k} \in X^k,\qquad d(\mathbf{u}) = \big( d(u_i, u_j) \big)_{0 \le i,j \le k}, \qquad
    u_0 = \rho.
\]
In the terminology of \cite{greven2009convergence, Depperschmidt2011},
$\Phi$ is called the \emph{monomial} corresponding to $\phi$. The space
of (equivalence classes of) mmm-spaces can be endowed with a topology
such that $(\mathcal{X}_n)_{n \ge 1}$ converges to $\mathcal{X}$ if and
only if 
\[
    \Phi(\mathcal{X}_n) \goesto{n \to \infty} \Phi(\mathcal{X})
\]
for all $\Phi$ of the form \eqref{eq:polynomial} with $\phi$ bounded and
continuous. This topology is called the (pointed marked)
\emph{Gromov-weak topology}. 

\paragraph{Gromov-vague topology.}
The Gromov-vague topology is obtained by restricting the Gromov-weak
topology to balls of finite radius around the pointed element of the
mmm-space. For an mmm-space $\mathcal{X}$ and $R > 0$, let
$\mathcal{X}^{(R)}$ be the restriction of $\mathcal{X}$ to the closed
ball $B_R(\mathcal{X})$ of $(X,d)$ with center $\rho$ and radius $R$.
Namely, $\mathcal{X}^{(R)}$ is the equivalence class of $(X, d, \rho,
\mu^{(R)})$, where
\[
    \forall A \in \mathscr{B}(X \times E),\quad 
    \mu^{(R)}(A) =  \mu\big( A \cap (B_R(\mathcal{X}) \times E) \big).
\]
Following the definition of \cite{athreya16}, a sequence
$(\mathcal{X}_n)_{n \ge 1}$ of mmm-spaces converges to $\mathcal{X}$ in
the (pointed marked) \emph{Gromov-vague topology} if
\[
    \mathcal{X}^{(R)}_n \goesto{n \to \infty} \mathcal{X}^{(R)}
\]
in the Gromov-weak topology for a.e.\ $R > 0$.
It will be convenient to reformulate the Gromov-vague convergence in
terms of functionals $\Phi$ of the form \eqref{eq:polynomial}.
We say that a function $\phi \colon \mathscr{D}^*_k \to \R$ has \emph{bounded
support} if there exists $R > 0$ such that
\begin{equation} \label{eq:boundedDistanceMat}
    \supp \phi \subseteq \mathscr{D}^*_{k,R} \coloneqq 
    \big\{ ((d_{ij})_{i,j}, (e_i)_i) \in \mathscr{D}^*_k : 
        \forall i \le k,\, d_{0i} \le R \big\}.
\end{equation}
The set appearing in the right-hand side is the set of distance matrices
such that the distance of each vertex to the pointed element (with index
$0$) is at most $R$. This definition is designed in such a way that, if
$\phi$ has bounded support and $\mathcal{X}$ is an mmm-space,
$\phi(d(\mathbf{u}), \mathbf{e}) = 0$ whenever $d(\rho, u_i) > R$ for
some $i \le k$. Let $\Pi_{\mathrm{b}}$ denote the set of all functionals
$\Phi$ of the form \eqref{eq:polynomial} with $k \ge 1$, such that $\phi$
is continuous, bounded, and with bounded support. It should be clear that
$\Phi(\mathcal{X}) = \Phi(\mathcal{X}^{(R)})$, and therefore
$(\mathcal{X}_n)_{n \ge 1}$ converges Gromov-vaguely to $\mathcal{X}$ if
and only if
\[
    \forall \Phi \in \Pi_{\mathrm{b}},\quad 
    \Phi(\mathcal{X}_n) \goesto{n \to \infty} \Phi(\mathcal{X}).
\]

\paragraph{Gromov--Hausdorff-weak topology.}
Finally, we define a straightforward extension to marked spaces of the
Gromov--Hausdorff-weak topology of \cite{athreya16}. We say that
$(\mathcal{X}_n)_{n \ge 1}$ converges to $\mathcal{X}_\infty$ in the
(pointed marked) \emph{Gromov--Hausdorff-weak topology} if there exist
\begin{enumerate}
    \item a complete separable metric space $(Z, d_Z)$;
    \item isometries $\iota_n \colon \{\rho_n\} \cup \supp_{X_n} \mu_n \to Z$,
        $n = 1, 2,\dots, \infty$;
\end{enumerate}
such that 
\[
    \iota_n(\rho_n) \goesto{n \to \infty} \iota_\infty(\rho_\infty),\qquad
    \iota_n(\supp_{X_n} \mu_n) \goesto{n \to \infty}
    \iota_\infty(\supp_{X_\infty} \mu_\infty),
\]
in $(Z,d_Z)$ and in the Hausdorff topology respectively, and 
\[
    \mu_n \circ \tilde{\iota}_n^{-1} \goesto{n \to \infty}
    \mu_\infty \circ \tilde{\iota}_\infty^{-1} 
\]
weakly as measures on $Z \times E$, where $\tilde{\iota}_n(x,e) =
(\iota_n(x), e)$. If the mmm-spaces have full support, this notion of
convergence coincides with the more standard convergence in the (pointed
marked) Gromov--Hausdorff--Prohorov topology \cite{miermont2009,
Abraham2013}. In this case, we can use the two notions of convergence
interchangeably, see the discussion in Section~5 of \cite{athreya16}. 
Although less standard, the Gromov--Hausdorff-weak topology is defined on
the same equivalence classes as the Gromov-weak topology, making it a
more natural choice when trying to reinforce a Gromov-weak convergence.

\subsection{Vague convergence and method of moments}

Each of the previous three topologies comes with a notion of weak
convergence of finite measures on $\mathscr{X}$. However, since the
limiting measure in our application (the free Brownian CRT) is infinite,
we need to work with an appropriate notion of \emph{vague convergence}.
We follow that considered in \cite{foutel2024vague}, which is obtained by
restricting the measures to sets of mmm-spaces with a large mass, of the
form $\{ \abs{X} \ge \eta \}$ for some $\eta > 0$. More formally, a
measure $\mathds{M}$ on $\mathscr{X}$ is called \emph{locally finite} if 
\[
    \forall \eta > 0,\quad \mathds{M}( \abs{X} \ge \eta ) < \infty.
\]
A sequence $(\mathds{M}_n)_{n \ge 1}$ of locally finite measures
converges \emph{vaguely} to $\mathds{M}$ for one of the above topologies
if 
\[
    \mathds{M}_n\big[ F(\mathcal{X}) ; \abs{X} \ge \eta \big] 
    \goesto{n \to \infty}
    \mathds{M}\big[ F(\mathcal{X}) ; \abs{X} \ge \eta \big],
\]
for any bounded function $F \colon \mathscr{X} \to \R$ which is
continuous in the corresponding topology. Here, we have used the
notation $\mu[F; A] = \int F(x) \indic_A(x) \mu(\diff x)$. We refer to
\cite{foutel2024convergence} for a more detailed exposition.

We can now state our method of moments, which provides a connexion
between convergence of the moments and vague convergence in the
Gromov-vague topology. It is a straightforward adaptation of the main
result in \cite{foutel2024convergence}.

\begin{proposition}[Method of moments] \label{prop:methodMoments}
    Let $(\mathds{M}_n)_{n \ge 1}$ be a sequence of locally finite
    measures on $\mathscr{X}$. Suppose that 
    \begin{equation} \label{eq:methodMoment}
        \forall \Phi \in \Pi_{\mathrm{b}},\quad
        \mathds{M}_n[\Phi] \goesto{n \to \infty} \mathds{M}[\Phi],
    \end{equation}
    for some limit measure $\mathds{M}$ on $\mathscr{X}$. Suppose in
    addition that 
    \begin{equation} \label{eq:carleman}
        \forall R > 0,\quad
        \sum_{k \ge 1} \mathds{M}\big[ \abs{\mathcal{X}^{(R)}}^k \big]^{-\frac{1}{2k}} = \infty.
    \end{equation}
    Then $(\mathds{M}_n)_{n \ge 1}$ converges to $\mathds{M}$ vaguely in
    the Gromov-vague topology.
\end{proposition}

\begin{proof}
    For $k \ge 1$ we can construct some measures $(\M^k_n)_{n \ge 1}$ and
    $\M^k$ on $\mathscr{D}^*_k$ such that
    \[
        \forall \phi \colon \mathscr{D}^*_k \to \R_+,\quad \M^k_n[\phi] = \mathds{M}_n[\Phi],
        \quad \M^k[\phi] = \mathds{M}[\Phi],
    \]
    where $\Phi$ is as in \eqref{eq:polynomial}. Recall the notation
    $\mathscr{D}^*_{k,R}$ in \eqref{eq:boundedDistanceMat}. Our
    assumption \eqref{eq:methodMoment} can be reformulated as saying that 
    \[
        \M^k_n \goesto{n \to \infty} \M^k
    \]
    in the sense of weak\textsuperscript{\#} convergence of measures of
    \cite[Section~A2.6]{daley2003introduction}, where a set $A \in \mathscr{D}^*_k$ is
    called bounded if $A \subseteq \mathscr{D}^*_{k,R}$ for some $R > 0$.
    (This convergence is also called vague convergence in
    \cite[Chapter~4]{kallenberg2017random}.) It is now standard (see 
    \cite[Proposition~A2.6.II]{daley2003introduction}) that the
    restriction of $\M^k_n$ to $\mathscr{D}^*_{k,R}$ converges
    \emph{weakly} to that of $\M^k$, for a.e.\ $R > 0$. Reformulating
    this in terms of mmm-spaces, we obtain that, for a.e.\ $R > 0$ and
    any continuous bounded $\phi \colon \mathscr{D}^*_k \to \R$ with
    corresponding polynomial $\Phi$,
    \[
        \mathds{M}_n\big[\Phi(\mathcal{X}^{(R)})\big]
        = \M^k_n\big[ \phi; \mathscr{D}^*_{k,R} \big]
        \goesto{n \to \infty}
        \M^k\big[ \phi; \mathscr{D}^*_{k,R} \big]
        =
        \mathds{M}\big[\Phi(\mathcal{X}^{(R)})\big].
    \]
    The method of moments for vague convergence in the Gromov-weak
    topology \cite[Theorem~4.2]{foutel2024vague} shows that the ``law''
    of $\mathcal{X}^{(R)}$ under $\mathds{M}_n$ converges vaguely to that
    under $\mathds{M}$. (Here we have used that Carleman's condition
    \eqref{eq:carleman} is fulfilled for each $\mathcal{X}^{(R)}$.) This
    proves the result.
\end{proof}

\begin{remark}
    An equivalent reformulation of the condition \eqref{eq:methodMoment}
    is that the measures $(\M^k_n)_{n \ge 1}$ appearing in the proof need
    to converge vaguely to $\M^k$. This is the point of view followed in
    \cite{foutel2024convergence}, and $\M^k$ is the $k$-th moment measure
    of $\mathds{M}$. These moment measures can be seen as the unplanar
    versions of the planar factorial moment defineds in \eqref{eq:defMoment}.
\end{remark}

\subsection{The Brownian CRT, CPP and their moments}
\label{SS:BrownianCRT}

In this section, we provide a brief construction of the Brownian CRT and
of the Brownian CPP to set some notation up and give an expression for
their moments. For our purpose, it is more natural to work with the
``free'' Brownian CRT (obtained under the infinite excursion measure)
than with the usual Brownian CRT obtained from a normalized excursion.

\paragraph{The Brownian CRT.}  We
refer to \cite{LeGall2005} for a careful exposition.
An excursion is a continuous function $f \colon \R_+ \to \R_+$ such that
there is some $\omega > 0$ (the length of the excursion) verifying that
$f(u) > 0$ for $u \in (0, \omega)$ and $f(u) = 0$ otherwise. Define a
(pseudo) metric on $[0, \omega]$ as
\[
    \forall u < v \in [0,\omega],\quad d_f(u,v) = d_f(v,u) 
    = f(u) + f(v) - 2 \inf_{w \in [u,v]} f(w).
\]
Up to taking a quotient, $\mathcal{T}_f \coloneqq ([0, \omega], d_f, 0, \Leb)$
is a metric measure space rooted at $0$. This construction defines a map
$C \colon f \mapsto \mathcal{T}_f$ (the contour map) from the space of
excursions to the space of metric measure spaces. If $a > 0$, let us also
define the rescaled contour map as $C_a(f) = C(af)$.

Following \cite{LeGall1993}, we let let $\mathbf{n}$ be Itô's excursion
measure of a variance one reflected Brownian motion, normalised in such
a way that 
\[
    \mathbf{n}\Big( \sup_{x \in [0,\omega]} f(x) \ge t \Big) = \frac{1}{2t}.
\]
We let the Brownian CRT with variance $\Sigma^2$ be the random variable
$\mathcal{T}_{\mathrm{b}} = (T_{\mathrm{b}}, d_{\mathrm{b}},
\rho_{\mathrm{b}}, \mu_{\mathrm{b}})$ with ``distribution'' 
\begin{equation} \label{eq:defCRT}
    \frac{2}{\Sigma} \big(\mathbf{n} \circ C_{2/\Sigma}^{-1}\big),
\end{equation}
which is obtained by taking the pushforward of the excursion measure $\mathbf{n}$
through the rescaling by $\Sigma/2$ of the contour map $C$. (Note that 
it corresponds to the excursion measure of a Brownian motion with
variance $4/\Sigma^2$.) In what follows, the expectation of some function
of $\mathcal{T}_{\mathrm{b}}$ should be interpreted as the integral of
that function against the (infinite) measure \eqref{eq:defCRT}.

For a probability measure $\pi$ on $E$, we let $\mathcal{T}_{\mathrm{b},\pi}$ 
be the mmm-space obtained by assigning independent marks distributed as
$\pi$ to the leaves of $\mathcal{T}_{\mathrm{b},\pi}$. More formally, it
is defined as 
\[
    \mathcal{T}_{\mathrm{b}, \pi} = \big( T_{\mathrm{b}}, d_{\mathrm{b}},
        \rho_{\mathrm{b}}, \mu_{\mathrm{b}} \otimes \pi \big).
\]
Before stating the result, recall that $\Lambda_k$ stands for the uniform
measure on binary trees with $k$ leaves defined in \eqref{eq:uniformTree}, and that
$D(\theta) \in \mathscr{D}_k$ is the $(k+1)\times (k+1)$ matrix of
pairwise distances between the $k$ leaves and the root of $\theta$
defined in \eqref{eq:defDistanceMatrix}. For a functional $\phi \colon
\mathscr{D}^*_k \to \R$ and a permutation $\sigma$ of $[k]$, we let
$\phi_\sigma$ be the functional obtained by relabeling the leaves
according to $\sigma$, namely
\[
    \forall \big((d_{ij})_{i,j}, (e_i)_i \big) \in \mathscr{D}^*_k,\quad
    \phi_\sigma\big( (d_{ij})_{i,j}, (e_i)_i \big) = 
    \phi\big((d_{\sigma_i, \sigma_j})_{i,j}, (e_{\sigma_i})_i\big), 
    \quad \sigma_0 = 0.
\]
The set of all permutations of $[k]$ is denoted by $\mathscr{P}_k$.
The following expression for the moment is merely a reformulation of
\cite[Theorem~3]{LeGall1993}.

\begin{proposition}[Moment of the Brownian CRT] \label{prop:momentMarkedCRT}
    For any functional $\phi \colon \mathscr{D}^*_k \to \R_+$,
    the moment of the Brownian CRT for the corresponding polynomial
    $\Phi$ is 
    \begin{equation*}
        \E\big[ \Phi( \mathcal{T}_{\mathrm{b},\pi} ) \big]
        = \Big( \frac{\Sigma^2}{2} \Big)^{k-1} \sum_{\sigma \in \mathscr{P}_k} 
        \int \phi_\sigma\big( D(\theta), (X_i)_{i\le k} \big) \Lambda_k(\diff \theta),
    \end{equation*}
    where $(X_i)_{i \le k}$ are i.i.d.\ and distributed as $\pi$.
\end{proposition}

\begin{proof}
    Since the sampling measure is of the form $\mu_{\mathrm{b}} \otimes
    \pi$, it is sufficient to prove the result for functionals $\phi$ that are
    independent of the marks, see \cite[Corollary~1]{foutel2023convergence}.
    For such a functional, if $\Phi$ is the corresponding polynomial in
    \eqref{eq:polynomial} and $\mathcal{T}_{\mathrm{b},\pi}$,
    \begin{align*}
        \E\big[ \Phi( \mathcal{T}_{\mathrm{b},\pi} ) \big]
        &= \frac{2}{\Sigma}
        \int \bigg(
            \int_{[0,\omega]^k} \phi\big( d_{2f/\Sigma}(\mathbf{u}) \big) \diff \mathbf{u} 
        \bigg) \diff \mathbf{n}(f) \\
        &= \frac{2}{\Sigma}
        \sum_{\sigma \in \mathscr{P}_k}
        \int \bigg(
        \int_{0\le u_1 < \dots < u_k \le \omega} \phi_\sigma\big( \tfrac{2}{\Sigma} d_f(\mathbf{u}) \big) \diff \mathbf{u} 
        \bigg) \diff \mathbf{n}(f).
    \end{align*}
    By \cite[Theorem~4 of Chapter~III]{legall1999spatial} (which is a reformulation of
    \cite[Theorem~3]{LeGall1993}), 
    \begin{align*}
        \frac{2}{\Sigma}
        \sum_{\sigma \in \mathscr{P}_k}
        \int \bigg(
        \int_{0\le u_1 < \dots < u_k \le \omega} 
        &\phi_\sigma\big( \tfrac{2}{\Sigma} d_f(\mathbf{u}) \big) \diff \mathbf{u} \bigg) \diff \mathbf{n}(f) \\
        &= 2^{k-1}
        \frac{2}{\Sigma} 
        \sum_{\sigma \in \mathscr{P}_k}
        \int \phi_\sigma\big( \tfrac{2}{\Sigma} D(\theta) \big) 
        \Lambda_k(\diff \theta) \\
        &= 
        \frac{2^k}{\Sigma} \Big(\frac{\Sigma}{2}\Big)^{2k-1}
        \sum_{\sigma \in \mathscr{P}_k}
        \int \phi_\sigma\big( D(\theta') \big) 
        \Lambda_k(\diff \theta'),
    \end{align*}
    where we have used the change of variable $\theta' = 2\theta /
    \Sigma$ in \eqref{eq:uniformTree}.
\end{proof}

\begin{remark}
    Note that the balls of the Brownian CRT fulfill Carleman's
    condition, since their moments can be bounded as 
    \[
        \E\big[ \abs{T^{(R)}_{\mathrm{b}}}^k \big] \le k!\, R^{2k-1} \Big( \frac{\Sigma^2}{2} \Big)^{k-1}.
    \]
\end{remark}

\paragraph{The Brownian CPP.} The Brownian coalescent point process
corresponds to the sphere at height $1$ of the Brownian CRT, under the
conditional measure that the CRT has height greater than $1$. By Itô's
theory of excursions, it can be constructed directly from a Poisson
point process $P$ on $\R_+ \times [0,1]$ with intensity $\diff u \otimes
\tfrac{\diff s}{s^2}$, which corresponds to the depths of the downwards
excursions below level $1$ of a Brownian motion, see \cite{Popovic2004,
foutel2023convergence} for more details. Define a distance on $\R_+$ as 
\[
    \forall u \le v,\quad d_P(u,v) = d_P(v,u) = 
    2 \sup \{ s : (w,s) \in P,\, u \le w \le v \}
\]
(The factor $2$ accounts for the fact that the graph distance is twice
the time to the most-recent common ancestor.) Let $Z$ be an independent
exponential random variable with mean $1$. The Brownian CPP with variance
$\Sigma^2$ is the random mmm-space $\mathcal{U}_{\mathrm{b}} =
(U_{\mathrm{b}}, d_{\mathrm{b}}, \rho_{\mathrm{b}}, \mu_{\mathrm{b}})$
which is the equivalence class in $\mathscr{X}$ of $([0,Z], d_P, \rho,
\tfrac{\Sigma^2}{2} \Leb)$, where $\rho$ is a point lying at distance $1$
from all other points in $[0,Z]$. 

As above, we let $\mathcal{U}_{\mathrm{b},\pi}$ be obtained by
assigning independent marks distributed as $\pi$ to the leaves of a
Brownian CPP, $\mathcal{U}_{\mathrm{b}, \pi} = \big(U_{\mathrm{b}}, 
d_{\mathrm{b}}, \rho_{\mathrm{b}}, \mu_{\mathrm{b}} \otimes \pi\big)$.
The expression below is proved in \cite[Corollary~1]{foutel2023convergence}, 
where we recall that $\widetilde{\Lambda}_k$ stands for the uniform
measure on ultrametric trees with $k$ leaves defined in
\eqref{eq:uniformUltrametric}.

\begin{proposition}[Moment of the Brownian CPP] \label{prop:momentCPP}
    For any functional $\phi \colon \mathscr{D}^*_k \to \R_+$,
    the moment of the Brownian CPP for the corresponding polynomial
    $\Phi$ is 
    \begin{equation*}
        \E\big[ \Phi(\mathcal{U}_{\mathrm{b}, \pi}) \big]
        = \Big( \frac{\Sigma^2}{2} \Big)^k \sum_{\sigma \in \mathscr{P}_k} 
        \int \E\big[ \phi_\sigma\big( D(\theta), (X_i)_{i \le k} \big) \big]
        \widetilde{\Lambda}_k(\diff \theta),
    \end{equation*}
    where $(X_i)_{i \ge 1}$ are i.i.d.\ and distributed as $\pi$.
\end{proposition}

\subsection{Tightness results}

Once convergence is established in the Gromov-vague topology (typically
through the method of moments) we will be interested in reinforcing it
to a Gromov--Hausdorff--Prohorov convergence. We provide two tightness
criteria that we need in this direction. The first one will be used to go
from a Gromov-vague to a Gromov-weak convergence, and amounts to
controlling the mass escaping at an infinite distance from the pointed
vertex. The second one is a simple adaptation of the tightness criterion
in the Gromov--Hausdorff--Prohorov topology of \cite{athreya16}. It relies on
controlling the \emph{lower mass function}, defined for $\delta > 0$ and
$\mathcal{X} \in \mathscr{X}$ as
\begin{equation} \label{eq:lowerMassFunc}
    m_{\delta}(\mathcal{X}) = \inf \{ \abs{B_{\delta,x}(X)} : (x,e) \in
    \supp \mu \},
\end{equation}
where $B_{\delta,x}(X)$ is the closed ball of radius $\delta$ and center
$x$ of $(X, d)$.

\begin{lemma}[Gromov-vague to Gromov-weak] \label{lem:Gv2Gw}
    Let $(\mathds{M}_n)_{n \ge 1}$ be a sequence of locally finite
    measures on $\mathscr{X}$ that converges vaguely to $\mathds{M}$ for
    the Gromov-vague topology. If
    \begin{equation} \label{eq:Gv2Gw1}
        \forall \epsilon > 0,\: \exists R > 0,\quad 
        \limsup_{n \ge 1} \mathds{M}_n\big( \abs{X \setminus B_R(X)} \ge
            \epsilon \big) \le \epsilon
    \end{equation}
    the convergence occurs vaguely in the Gromov-weak topology.
\end{lemma}

\begin{proof}
    The result follows from simple but tedious manipulations. Let
    $d_{\mathrm{GP}}$ be the (pointed marked) Gromov--Prohorov distance,
    which can defined by adapting \cite[Definition~3.1]{Depperschmidt2011} 
    to take the pointed element into account. It is sufficient for our
    purpose to know that it induces the Gromov-weak topology and that 
    \begin{equation} \label{eq:dGPBound}
        d_{\mathrm{GP}}\big( \mathcal{X}, \mathcal{X}^{(R)} \big)
        \le \abs{X \setminus B_R(X)},
    \end{equation}
    see for instance \cite[Lemma~3.3]{boenkost2022genealogy}. Fix a
    bounded $F \colon \mathscr{X} \to \R$ which is $L$-Lipschitz for the
    distance $1 \wedge d_{\mathrm{GP}}$. For $\eta,\eta' > 0$ such that
    $\eta' < \eta$,
    \begin{multline} \label{eq:Gv2GW1}
        \mathds{M}_n\big[ F(\mathcal{X}) ; \abs{X} > \eta \big]
        =
        \mathds{M}_n\big[ F(\mathcal{X}^{(R)}) ; \abs{X^{(R)}} > \eta' \big]  \\
        +
        \mathds{M}_n\big[ F(\mathcal{X})-F(\mathcal{X}^{(R)})
        ; \abs{X^{(R)}} > \eta' \big] 
        +
        \mathds{M}_n\big[ F(\mathcal{X}) ; \abs{X} > \eta, \abs{X^{(R)}} \le
        \eta' \big]. 
    \end{multline}
    By the Gromov-vague convergence
    \[
        \lim_{R \to \infty} 
        \adjustlimits
        \lim_{\eta' \uparrow \eta} 
        \lim_{n \to \infty}
        \mathds{M}_n\big[ F(\mathcal{X}^{(R)}); \abs{X^{(R)}} > \eta' \big] 
        =
        \mathds{M}\big[ F(\mathcal{X}) ; \abs{X} \ge \eta \big],
    \]
    and the result is proved if we can show that the remaining two terms
    in \eqref{eq:Gv2GW1} vanish when limits are taken in the same order.
    For the last term,
    \[
        \mathds{M}_n\big[ F(\mathcal{X}) ; \abs{X} > \eta, \abs{X^{(R)}} \le \eta' 
        \big] \le 
        \norm{F}_\infty \mathds{M}_n\big( \abs{X \setminus B_R(X)} \ge \eta - \eta' \big),
    \]
    which by \eqref{eq:Gv2Gw1} vanishes under our limits. For the second term, for any $\epsilon > 0$,
    \begin{align*}
        \abs[\big]{ \mathds{M}_n\big[ 
            F(\mathcal{X})-F(\mathcal{X}^{(R)})
            ; \abs{X^{(R)}} > \eta'
        \big] }
        &\le L \mathds{M}_n\big[ 
            1 \wedge d_{\mathrm{GP}}\big(\mathcal{X}, \mathcal{X}^{(R)}\big) 
            ; \abs{X^{(R)}} > \eta' 
        \big]
        \\
        &\le L \mathds{M}_n\big[ 
            1 \wedge \abs{X \setminus B_R(X)} 
            ; \abs{X^{(R)}} > \eta' 
        \big] \\
        &\le \epsilon L \mathds{M}_n\big( \abs{B_R(X)} > \eta' \big)
        +  L \mathds{M}_n\big( \abs{X \setminus B_R(X)} > \epsilon \big),
    \end{align*}
    where we have used \eqref{eq:dGPBound} in the second line. We can end
    the proof by noting that, again by \eqref{eq:Gv2Gw1},
    \[
        \limsup_{R \to \infty} \limsup_{n \to \infty} \mathds{M}_n\big(
        \abs{X \setminus B_R(X)} > \epsilon \big) \le \epsilon
    \]
    and that 
    \[
        \limsup_{R \to \infty} \limsup_{\eta' \uparrow \eta} 
        \limsup_{n \to \infty} \mathds{M}_n\big( \abs{B_R(X)} > \eta' \big) 
        \le \mathds{M}\big( \abs{X} \ge \eta \big)
    \]
    and letting $\epsilon \to 0$.
\end{proof}

\begin{lemma}[Gromov-weak to Gromov--Hausdorff-weak] \label{lem:Gw2GHP}
    Suppose that $(\mathds{M}_n)_{n \ge 1}$ is a sequence of locally
    finite measures converging to $\mathds{M}$ vaguely in the Gromov-weak
    topology. Suppose that
    \begin{equation} \label{eq:Gw2GHw}
        \forall \eta, \delta > 0,\quad 
        \adjustlimits
        \lim_{\eta' \to 0} \sup_{n \ge 1} 
        \mathds{M}_n\big( \abs{X} > \eta, m_\delta(\mathcal{X}) < \eta' \big) = 0
    \end{equation}
    then $(\mathds{M}_n)_{n \ge 1}$ also converges vaguely in the
    Gromov--Hausdorff-weak topology.
\end{lemma}

\begin{proof}
    For $\epsilon > 0$ and $k \ge 1$, let $\epsilon_k = 2^{-k} \epsilon$
    and $\delta_k > 0$ be such that $\delta_k \to 0$. By our assumption,
    there exists $\eta'_k > 0$ such that 
    \[
        \sup_{n \ge 1} 
        \mathds{M}_n\big( \abs{X} > \eta, m_{\delta_k}(\mathcal{X}) <
        \eta_k' \big) \le \epsilon_k.
    \]
    We deduce that 
    \[
        \sup_{n \ge 1} 
        \mathds{M}_n\big( \{\abs{X} > \eta\} \setminus A \big)
        \le \epsilon, \qquad 
        A \coloneqq \bigcap_{k \ge 1} \{ m_{\delta_k}(\mathcal{X}) \ge \eta_k' \}.
    \]
    Using the vague convergence and, for instance, 
    \cite[Theorem~4.2]{kallenberg2017random} we can find a compact set
    $\mathscr{K} \subset \mathscr{X}$ for the Gromov-weak topology such
    that 
    \[
        \sup_{n \ge 1} 
        \mathds{M}_n\big( \{\abs{X} > \eta\} \setminus \mathscr{K} \big)
        \le \epsilon.
    \]
    Therefore,
    \[
        \sup_{n \ge 1} 
        \mathds{M}_n\big( \{\abs{X} > \eta\} \setminus (\mathscr{K} \cap A)\big)
        \le 2 \epsilon,
    \]
    which proves tightness of $(\mathds{M}_n)_{n \ge 1}$ (and hence
    the result) in the Gromov--Hausdorff-weak topology provided that we
    show that $\mathscr{K} \cap A$ is compact in that topology. Let
    $(\mathcal{X}_p)_{p \ge 1} \in \mathscr{K} \cap A$. Since
    $\mathscr{K}$ is compact we can assume that it admits a subsequence
    that converges to some limit $\mathcal{X}$ in the Gromov-weak topology.
    For each $\delta > 0$, $m_{\delta}(\mathcal{X}_p) > \eta'_k$ for $k$
    sufficiently large. Hence, \cite[Theorem~6.1]{athreya16} shows that
    the convergence of the subsequence also holds Gromov--Hausdorff-weakly.
\end{proof}

\begin{remark}
    If, for each $n \ge 1$, $\mathds{M}_n$ puts no mass on non-compact
    mmm-spaces, \cite[Lemma~3.4]{athreya16} and the dominated convergence
    theorem entail that 
    \[
        \forall \eta, \delta > 0,\quad 
        \lim_{\eta' \to 0} \mathds{M}_n\big( \abs{X} > \eta, 
        m_\delta(\mathcal{X}) < \eta' \big) = 0
    \]
    and thus \eqref{eq:Gw2GHw} can be replaced by the simpler condition
    \begin{equation*}
        \forall \eta, \delta > 0,\quad 
        \adjustlimits
        \lim_{\eta' \to 0} \limsup_{n \ge 1} 
        \mathds{M}_n\big( \abs{X} > \eta, m_\delta(\mathcal{X}) < \eta' \big) = 0.
    \end{equation*}
\end{remark}

\section{Proof of the invariance principles}
\label{S:proof}

We prove the invariance principles by following the approach that was
outlined in the introduction. First, in Section~\ref{SS:momentComputation}, 
we compute the moments of the branching process using the many-to-few
formula. We deduce convergence in the Gromov-vague topology
(Theorem~\ref{thm:GromovVague}) from this computation and the method of
moments in Section~\ref{SS:proofGvague}. Section~\ref{SS:tightness}
contains results of a technical nature which will show tightness of the
sequence in the Gromov--Hausdorff-weak topology, which are used in
Section~\ref{SS:proofGHP} to deduce the Gromov--Hausdorff--Prohorov
convergence of Theorem~\ref{thm:GHP}.

\subsection{Moments asymptotic}
\label{SS:momentComputation}

We first compute the limit of the moments of the branching
process. Recall that we have defined the moment measure $\M^k_x$ of the
branching process in \eqref{eq:momentMeasure} as a finite measure on the
set $\mathscr{T}_k^*$ of discrete marked planar trees with $k$ leaves.
By viewing a discrete tree as a continuous tree, $\M^k_x$ can be seen as
a measure on $\mathscr{T}^*_{k,\mathrm{c}}$ and we let $\bar{\M}^k_{x,n}$ 
be obtained by rescaling mass and edge lengths as follows
\begin{equation} \label{eq:rescaledMoment}
    \bar{\M}^k_{x,n}[F] \coloneqq \frac{1}{n^{2k}} 
    \M^k_x\big[ F\big( \tfrac{1}{n}\theta^* \big) \big]
    =
    \frac{1}{n^{2k}}
    \E\bigg[ \sum_{\substack{\mathbf{v} \in T^k \\ v_1 < \dots < v_k}} 
    F\big( \tfrac{1}{n}\tau_{\mathbf{v}}, X_{\mathbf{v}} \big) \bigg],
\end{equation}
for $F \colon \mathscr{T}^*_{k,\mathrm{c}} \to \R_+$. Note that
$\bar{\M}^k_{x,n}$ corresponds to the moment measure of the rescaled
metric measure space $\bar{\mathcal{T}}_n$ defined in
\eqref{eq:rescalingMmmSpace}. Let us define the corresponding measure for
the ultrametric measure space $\bar{\mathcal{U}}_n$ of \eqref{eq:defUMS}
obtained by restricting the summation to vertices at height $n$, namely
\begin{equation*} 
    \widetilde{\M}^k_{x,n}[F] = \frac{1}{n^k}
    \E\bigg[ \sum_{\substack{\mathbf{v} \in U_n^k \\ v_1 < \dots < v_k}} 
    F\big( \tfrac{1}{n}\tau_{\mathbf{v}}, X_{\mathbf{v}} \big) \bigg].
\end{equation*}
The objective of this section is to derive the following asymptotic for
these moment measures, where we recall that $\Lambda_k$ stands for the
uniform measure on continuous binary trees defined in \eqref{eq:uniformTree}
and $\widetilde{\Lambda}_k$ for that on binary ultrametric trees with height
$1$, see \eqref{eq:uniformUltrametric}.

\begin{proposition} \label{prop:momentComputation}
Under Assumption~\ref{ass:criticality}, for any $R > 0$ and any
continuous bounded $F \colon \mathscr{T}^*_{k,\mathrm{c}} \to \R_+$ whose
support is included in the set of trees of height at most $R > 0$,
\begin{align*}
    \lim_{n \to \infty} n \bar{\M}^k_{x,n}[F] 
    &= h(x) \Big( \frac{\Sigma^2}{2} \Big)^{k-1}
    \int_{\mathscr{T}^*_{k,\mathrm{c}}}
    \E\big[ F\big(\theta, (X_i)_i \big) \big] \Lambda_k(\diff \theta), \\
    \lim_{n \to \infty} n \widetilde{\M}^k_{x,n}[F] 
    &= h(x) \Big( \frac{\Sigma^2}{2} \Big)^{k-1}
    \int_{\mathscr{T}^*_{k,\mathrm{c}}}
    \E\big[ F\big(\theta, (X_i)_i \big) \big] \widetilde{\Lambda}_k(\diff \theta),
\end{align*}
uniformly in $x \in E$, and where $(X_i)_{i \le k}$ are i.i.d.\ and
distributed as $\pi$.
\end{proposition}

We start by showing that the types at the branch points and leaves of the
tree-indexed Markov chain are asymptotically independent and distributed
as $h \pi$. In all this section, we use as biasing function $\psi \equiv
h$ and remove the superscript from all the quantities defined in
Section~\ref{S:moments}.

\begin{lemma} \label{lem:typeConvergence}
    Let $(\tau_n)_{n \ge 1}\in \mathscr{T}_k$ be a sequence of discrete
    trees with $k$ leaves such that
    \[
        \tfrac{1}{n} \tau_n \goesto{n \to \infty} \theta
    \]
    for some binary continuous tree $\theta \in \mathscr{T}_{k,\mathrm{c}}$.
    Then for any continuous bounded $F \colon \mathscr{T}^*_{k,\mathrm{c}} \to \R$
    \[
        h(x) \Q_{x, \tau_n}\big[ \Delta_k F \big] \goesto{n \to \infty}
        h(x) \Big( \frac{\Sigma^2}{2} \Big)^{k-1} \E\big[ F(\theta, (X_i)_{i \le k}) \big]
    \]
    uniformly in $x \in E$, where $(X_i)_{i \ge 1}$ are i.i.d.\ and
    distributed as $\pi$.
\end{lemma}

\begin{proof}
Since the trees are converging deterministically, it is sufficient to
prove the result for functionals that only depend on the types of the
leaves, namely
\[
    \forall \theta^* = (\theta, (x_i)_{i \le k}),\quad F(\theta^*) = \prod_{i=1}^k f_i(x_i),
\]
for some continuous bounded $f_i \colon E \to \R$, $i \le k$.
We prove the result by an induction on the number of leaves $k$. For $k =
1$, the tree $\tau_n$ is made of a single branch of length $\ell^n_1$
with $\ell^n_1 \to \infty$. The many-to-one formula
(Proposition~\ref{prop:manyToOne}) applied to the bias $\psi \equiv h$
yields that
\begin{equation*} 
    h(x) \Q_{x,\tau_n}[\Delta_1 F] 
    = h(x) \E_x\Big[ \frac{f_1(\zeta_{\ell^1_n})}{h(\zeta_{\ell^1_n})} \Big]
    = \E_x\Big[ \sum_{v \in T, \abs{v} = \ell^1_n} f_1(X_v) \Big],
\end{equation*}
and the result is the content of our assumption \eqref{eq:convMean}. 

Suppose that the result holds for some $k \ge 1$ and let $v^n_1 < \dots <
v^n_{k+1}$ be the leaves of $\tau_n$. Since $\theta$ is binary, $\tau_n$ is
also binary for $n$ large enough. In that case, the tree spanned by
$(v^n_1, \dots, v^n_k, v^n_k \wedge v^n_{k+1})$ -- obtained by removing
the two external branches leading to $v^n_k$ and $v^n_{k+1}$ -- has $k$
leaves. Let us denote this subtree by $\widetilde{\tau}_n$, and let
$\widetilde{\tau}^*_n = (\widetilde{\tau}_n, (X_u)_{u \in
\widetilde{\tau}_n})$ be the corresponding marked tree. By construction
of the tree-indexed Markov chain we know that, conditional on
$\widetilde{\tau}^*_n$,
\[
    (X_{v^n_k}, X_{v^n_{k+1}})
    \overset{\mathrm{(d)}}{=}
    \big(
        \zeta_{\ell^n_k-b^n_k-1}, \zeta'_{\ell^n_{k+1}-b^n_k-1}
    \big),
\]
where $(\zeta_n)_{n \ge 0}$ and $(\zeta'_n)_{n \ge 0}$ are
two copies of the spine process, with starting points distributed as 
$(\chi_{1,2}, \chi_{2,2})$ under $\P_{X_{v^n_k \wedge v^n_{k+1}}}$
defined in \eqref{eq:jumpDistr}, and run independently.

Unfolding the definitions, 
\begin{multline*}
    \Q_{x,\tau_n}
    \Big[ \frac{m_2(X_{v^n_k \wedge v^n_{k+1}}) f_k(X_{v^n_k})f_{k+1}(X_{v^n_{k+1}}) }{h(X_{v^n_k})h(X_{v^n_{k+1}})}
    \:\Big|\: \widetilde{\tau}^*_n \Big] \\
    = 
    \E_{X_{v^n_k \wedge v^n_{k+1}}}\bigg[
        \sum_{i,j=1, i\ne j}^{\abs{\Xi}} 
        h(\xi_i) \E_{\xi_i}\Big[ \frac{f_k(\zeta_{\ell^n_k-b^n_k-1})}{h(\zeta_{\ell^n_k-b^n_k-1})} \Big]
        h(\xi_j) \E_{\xi_j}\Big[ \frac{f_{k+1}(\zeta_{\ell^n_{k+1}-b^n_k-1})}{h(\zeta_{\ell^n_{k+1}-b^n_k-1})} \Big]
    \bigg].
\end{multline*}
By \eqref{eq:convMean}, 
\[
    h(y) \E_{y}\Big[ \frac{f_k(\zeta_{\ell^n_k-b^n_k-1})}{h(\zeta_{\ell^n_k-b^n_k-1})} \Big]
    \goesto{n \to \infty}
    h(y) \angle{\pi, f_k}
\]
uniformly in $y \in E$. Moreover, since
\[
    \sup_{y \in E} \E_y[ \abs{\Xi}^2 ] < \infty,
\]
we deduce that
\begin{multline*}
    \E_y\bigg[
        \sum_{i,j=1, i\ne j}^{\abs{\Xi}} 
        h(\xi_i) \E_{\xi_i}\Big[ \frac{f_k(\zeta_{\ell^n_k-b^n_k-1})}{h(\zeta_{\ell^n_k-b^n_k-1})} \Big]
        h(\xi_j) \E_{\xi_j}\Big[ \frac{f_{k+1}(\zeta_{\ell^n_{k+1}-b^n_k-1})}{h(\zeta_{\ell^n_{k+1}-b^n_k-1})} \Big]
    \bigg] \\
    \goesto{n \to \infty}
    \angle{\pi, f_k} \angle{\pi, f_{k+1}}
    \E_y\Big[ \sum_{i,j=1, i\ne j}^{\abs{\Xi}} h(\xi_i) h(\xi_j) \Big]
    = 
    \angle{\pi, f_k} \angle{\pi, f_{k+1}} m_2(y)
\end{multline*}
uniformly in $y \in E$. By decomposing $\Delta_{k+1}$ as 
\[
    \Delta_{k+1}(\tau_n^*) = \Delta_k(\widetilde{\tau}_n^*)
    \frac{m_2(X_{v^n_k \wedge v^n_{k+1}}) }{2h(X_{v^n_k})h(X_{v^n_{k+1}})},
\]
our induction and the above uniform convergence show that
\begin{align*}
    h(x) &\Q_{x,\tau_n}[\Delta_{k+1} F]  \\
    &= h(x) \Q_{x,\tau_n}\bigg[ \Delta_k(\widetilde{\tau}^*_n)
        \prod_{i=1}^{k-1} f_i(X_{v^n_i}) 
        \Q_{x,\tau_n}\Big[  
            \frac{m_2(X_{v^n_k \wedge v^n_{k+1}})f_k(X_{v^n_k})f_{k+1}(X_{v^n_{k+1}}) }{2h(X_{v^n_k})h(X_{v^n_{k+1}})}
        \:\Big|\: \widetilde{\tau}^*_n \Big]
    \bigg] \\
    &\goesto{n \to \infty} h(x) \Big( \frac{\Sigma^2}{2} \Big)^{k-1}
    \prod_{i=1}^{k-1} \angle{\pi, f_i}
    \times 
    \frac{\angle{\pi, m_2}}{2} \angle{\pi, f_k} \angle{\pi, f_{k+1}} \\
    &= h(x) \Big( \frac{\Sigma^2}{2} \Big)^k
    \prod_{i=1}^{k+1} \angle{\pi, f_i},
\end{align*}
uniformly for $x \in E$.
\end{proof}

\begin{proof}[Proof of Proposition~\ref{prop:momentComputation}]
Let us proof the first part of the statement. By the many-to-few formula
(Theorem~\ref{thm:manyToFew}) applied to $\psi \equiv h$, for any
functional $F \colon \mathscr{T}^*_{k,c} \to \R_+$ we have that
\begin{align}
    n \bar{\M}^k_{x,n}[F]
    &= \frac{1}{n^{2k-1}} \sum_{\tau \in \mathscr{T}_k} h(x)
    \Q_{x,\tau}\big[\Delta_k \cdot F\big( \tfrac{1}{n} \theta^* \big) \big] \nonumber\\
    &= \frac{1}{n^{2k-1}} \sum_{\bm{\ell} \in \N^k} \sum_{\mathbf{b} \in \N^{k-1}}
    \prod_{i=1}^{k-1} \indic_{\{ b_i < \ell_i \wedge \ell_{i+1}\}}
    h(x) \Q_{x,\tau(\bm{\uell}, \mathbf{b})}\big[\Delta_k \cdot F\big( \tfrac{1}{n}\theta^* \big) \big] \nonumber\\ 
    &= \int_{\R_+^k}  \int_{\R_+^{k-1}} 
    \prod_{i=1}^{k-1} \indic_{\{ b_i < \ell_i \wedge \ell_{i+1}\}}
    h(x) \Q_{x,\tau_n(\bm{\uell}, \mathbf{b})}\big[ \Delta_k \cdot F\big(\tfrac{1}{n} \theta^*\big) \big]
    \,\diff \bm{\uell}\, \diff \mathbf{b}, \label{eq:computingMoment1}
\end{align}
where in the second line it is understood that $\tau(\bm{\uell},
\mathbf{b})$ is the tree corresponding to the heights $(\bm{\uell},
\mathbf{b})$ and in the last line we let $\tau_n(\bm{\uell}, \mathbf{b})$
be the (discrete) tree with leaf heights $\floor{n \bm{\uell}}$ and
branch point heights $\floor{n \mathbf{b}}$. 

Let $\theta(\bm{\uell}, \mathbf{b})$ be the (continuous) tree
with leaf and branch point heights $(\bm{\uell}, \mathbf{b})$. Clearly,
\[
    \frac{1}{n} \tau_n(\bm{\uell}, \mathbf{b}) 
    \goesto{n \to \infty}
    \theta(\bm{\uell}, \mathbf{b})
\]
and $\theta(\bm{\uell}, \mathbf{b})$ is a binary tree for a.e.\
$(\bm{\uell}, \mathbf{b})$. Lemma~\ref{lem:typeConvergence} shows that 
\[
    h(x) \Q_{x,\tau_n(\bm{\uell}, \mathbf{b})}\big[ \Delta_k \cdot F\big(\tfrac{1}{n} \theta^*\big) \big]
    \goesto{n \to \infty}
    h(x) \Big( \frac{\Sigma^2}{2} \Big)^{k-1} \E\big[ F(\theta, (X_i)_{i \le k}) \big].
\]
Therefore, our result will follow by dominated convergence provided 
the integrand in \eqref{eq:computingMoment1} is bounded uniformly in
$(\bm{\uell}, \mathbf{b})$. (Recall that the support of $F$ is included
in the set of trees with height at most $R$, so that the integral in
\eqref{eq:computingMoment1} can be restricted to a bounded domain.)

It is sufficient to prove by an induction on $k$ that 
\begin{equation} \label{eq:computingMoment2}
    \sup_{\tau \in \mathscr{T}_k} \sup_{x \in E} h(x)
    \Q_{x,\tau}\big[ \Delta_k \big] < \infty.
\end{equation}
For $k = 1$, if $\tau$ is the tree made of a single branch of length
$\ell+1$, 
\[
    h(x) \Q_{x,\tau}\big[ \Delta_1 \big] 
    = \E_x[ Z_\ell ]
    \goesto{\ell \to \infty} h(x),
\]
uniformly in $x \in E$ by \eqref{eq:convMean}, where $Z_\ell$ is the number of
particles at time $\ell$. Therefore \eqref{eq:computingMoment2} holds for
$\ell$ large enough. For $\ell$ fixed, we use that 
\[
    h(x) \Q_{x,\tau}\big[\Delta_1 \big] = \E_x[Z_\ell] = \E_x\Big[
    \sum_{u \in T, \abs{u} = \ell-1} \E_{X_u}[\abs{\Xi}] \Big]
    \le \big( \sup_{y \in E} \E_y[\abs{\Xi}] \big) \E_x[ Z_{\ell-1}]
\]
and an induction on $\ell$ to obtain \eqref{eq:computingMoment2}.

For a general $k \ge 1$ and $\tau \in \mathscr{T}_k$, note that there
exists $C_k$ independent of $y$ and $\tau$ such that
\[
    m_{\abs{\pi(\tau)}}(y)
    \E_{y}\Big[ \prod_{i=1}^{\abs{\pi(\tau)}}
    \Q_{\chi_i, S_i(\tau)}[ \Delta_{\abs{\pi_i}} ]
    \Big]
    =
    \E_y\Big[ \sum_{\substack{i_1, \dots, i_d = 1\\\text{$(i_p)$ distinct}}}^{\abs{\Xi}}
    \prod_{i=1}^d h(\xi_i) \Q_{\xi_i, S_i(\tau)}[ \Delta_{\abs{\pi_i}} ]
    \Big] < C_k,
\]
where we have used our induction and the moment bound
\eqref{eq:momentBound}. Finally, by construction of the tree-indexed
Markov chain, 
if $w$ is the first branch point of $\tau$ and $\abs{w} = p$,
\begin{align*}
    h(x) \Q_{x,\tau}[\Delta_k]
    &= h(x) \E_x\Big[
        \frac{m_{\abs{\pi(\tau)}}(\zeta_p)}{\abs{\pi(\tau)}!\, h(\zeta_p)}
        \E_{\zeta_p}\Big[ \prod_{i=1}^{\abs{\pi(\tau)}}
        \Q_{\chi_i, S_i(\tau)}[ \Delta_{\abs{\pi_i}} ] \Big]
    \Big] \\
    &\le C_k h(x) \E_x\Big[ \frac{1}{h(\zeta_p)} \Big],
\end{align*}
and we conclude by using the case $k=1$.

We now prove the second part of the statement on ultrametric spaces.
By Lemma~\ref{lem:cppConstruction} or the discussion at the end of
Section~\ref{SS:constructionTree}, the set of all ultrametric trees with
$k$ leaves at height $n$ is simply obtained by choosing the depth of the
$k-1$ successive branching points in $\{0, \dots, n-1\}$. Applying the
many-to-few formula and discarding all trees which are not ultrametric,
we obtain that 
\begin{align*}
    n \widetilde{\M}^k_{x,n}[F]
    &= \frac{1}{n^{k-1}} \sum_{b_1,\dots,b_{k-1}=0}^{n-1}
    h(x) \Q_{x,\tau(\mathbf{b})}\big[\Delta_k \cdot F\big( \tfrac{1}{n}\theta^* \big) \big] \nonumber\\ 
    &= \int_{[0,1]^{k-1}}
    h(x) \Q_{x,\tau_n(\mathbf{b})}\big[ \Delta_k \cdot F\big(\tfrac{1}{n} \theta^*\big) \big]
    \diff \mathbf{b},
\end{align*}
where $\tau(\mathbf{b})$ and $\tau_n(\mathbf{b})$ are the ultrametric tree
with $k$ leaves at height $n$ whose branch points are at height  
$\mathbf{b}$ and $\floor{n\mathbf{b}}$ respectively. The result
follows by dominated convergence as above, using Lemma~\ref{lem:typeConvergence} 
and the bound \eqref{eq:computingMoment2}.
\end{proof}

\subsection{Proof of Theorem~\ref{thm:GromovVague}}
\label{SS:proofGvague}

We are ready to prove Theorem~\ref{thm:GromovVague}. For later purpose,
we will need the convergences in that result to hold uniformly in the
initial condition $x \in E$. By this, we mean that the convergence holds
uniformly with respect to any metric that induces the vague convergence on the
Gromov-vague or Gromov-weak topology. For examples of such metrics, see
\cite[Equation~(5.1)]{foutel2024vague} for the Gromov-weak topology, or
\cite[Equation~(A2.6.1)]{daley2003introduction} for general vague convergence.

\begin{proposition} \label{prop:firstStep}
    Suppose that  Assumption~\ref{ass:criticality} holds. Then
    \begin{equation} \label{eq:vagueLimits}
        n \mathscr{L}_x( \bar{\mathcal{T}}_n )
        \goesto{n \to \infty}
        h(x) \mathscr{L}( \mathcal{T}_{\mathrm{b}, \pi} ),
        \quad n \mathscr{L}_x( \bar{\mathcal{U}}_n )
        \goesto{n \to \infty}
        \frac{2h(x)}{\Sigma^2} \mathscr{L}( \mathcal{U}_{\mathrm{b}, \pi} ),
    \end{equation}
    uniformly in $x \in E$, vaguely in the Gromov-vague and Gromov-weak
    topology respectively.
\end{proposition}

\begin{proof}
    Let $\phi \colon \mathscr{D}^*_k \to \R$ be continuous bounded with
    support included in $\mathscr{D}^*_{k,R}$ (defined in
    \eqref{eq:boundedDistanceMat}) for some $R > 0$. Fix a sequence
    $(x_n)_{n \ge 1} \in E$. Since $h$ is bounded, up to extracting a
    subsequence, we assume that $h(x_n) \to h(x)$ for some $x \in E$. 

    We start with the convergence of $\bar{\mathcal{T}}_n$. The moment of
    $\bar{\mathcal{T}}_n$ evaluated at the polynomial
    $\Phi$ corresponding to $\phi$ is
    \begin{equation} \label{eq:proofVague1}
        \E_{x_n}\big[ \Phi(\bar{\mathcal{T}}_n) \big] 
        = \frac{1}{n^{2k}} \E_{x_n}\Big[ 
            \sum_{\mathbf{v} \in T^k} 
            \phi\big( \tfrac{d(\mathbf{v})}{n}, X_{\mathbf{v}} \big)
        \Big].
    \end{equation}
    We will compare this expression to the rescaled moment measure
    $\bar{\M}^k_{{x_n},n}$ in \eqref{eq:rescaledMoment} by
    removing from the sum all terms that do not span a tree with exactly
    $k$ leaves. For a tree $\tau$, let $A_k(\tau)$ be the number of
    $k$-tuples of vertices $\mathbf{v} = (v_1, \dots, v_k) \in \tau$ such that
    $\tau_\mathbf{v}$ has strictly less than $k$ leaves. By removing
    these terms from \eqref{eq:proofVague1} and re-ordering the vertices
    in increasing planar order,
    \begin{equation} \label{eq:firstStep1}
        \E_{x_n}\big[ \Phi(\bar{\mathcal{T}}_n) \big] 
        = \sum_{\sigma \in \mathscr{P}_k}
        \bar{\M}^k_{{x_n},n}[ \phi_\sigma \circ D]
        + B_{{x_n},n},
    \end{equation}
    where $B_{{x_n},n}$ can be bounded as 
    \[
        \abs{B_{{x_n},n}} \le \frac{\norm{\phi}_\infty}{n^{2k}} \E_{x_n}[ A_k(T^{(Rn)}) ].
    \]
    By Lemma~\ref{lem:factorial2Moment}, since $T^{(Rn)}$ has height at
    most $Rn$ we obtain the further bound
    \[
        \abs{B_{{x_n},n}} \le \frac{\norm{\phi}_\infty}{n^{2k}} C_k Rn
        \E_{x_n}[ \abs{T^{(Rn)}}^{k-1} ]
        = \frac{C'}{n^2} \E_{x_n}\big[ \abs{\bar{T}_n^{(R)}}^{k-1} \big].
    \]
    It will follow by an induction on $k$ that
    \begin{equation}\label{eq:inductionFactorial}
        \sup_{n \ge 1} \sup_{x \in E}
        n \E_x\big[ \abs{\bar{T}_n^{(R)}}^{k-1} \big] 
        < \infty.
    \end{equation}
    Therefore, by Proposition~\ref{prop:momentComputation}, letting $n
    \to \infty$ in \eqref{eq:firstStep1} shows 
    \[
        n \E_{x_n}\big[ \Phi(\bar{\mathcal{T}}_n) \big] 
        \goesto{n \to \infty}
        h(x) \Big( \frac{\Sigma^2}{2} \Big)^{k-1}
        \sum_{\sigma \in \mathscr{P}_k}
        \int 
        \E\big[ \phi_\sigma\big( D(\theta), (X_i)_i \big) \big] \Lambda_k(\diff
        \theta),
    \]
    where $(X_i)_{i \le k}$ are i.i.d.\ and distributed as $\pi$.
    The right-hand side is the moment of the marked Brownian CRT with
    variance $\Sigma^2$ computed in Proposition~\ref{prop:momentMarkedCRT}.
    The method of moments (Proposition~\ref{prop:momentComputation}) shows
    that 
    \[
        n \mathscr{L}_{x_n}( \bar{\mathcal{T}}_n )
        \goesto{n \to \infty}
        h(x) \mathscr{L}( \mathcal{T}_{\mathrm{b},\pi} )
    \]
    vaguely in the Gromov-vague topology. (It is straightforward to check
    from the explicit expression of the moments that Carleman's condition
    \eqref{eq:carleman} is fulfilled.) This also shows that
    \eqref{eq:inductionFactorial} holds for $k+1$ and the induction is
    complete.

    The convergence of $\bar{\mathcal{U}}_n$ follows along the same
    lines. Re-ordering the terms in the sum,
    \begin{align*}
        \E_{x_n}\big[ \Phi(\bar{\mathcal{U}}_n) \big] 
        &= \frac{1}{n^k} \E_{x_n}\Big[ 
            \sum_{\mathbf{v} \in U_n^k} 
            \phi\big( \tfrac{d(\mathbf{v})}{n}, X_{\mathbf{v}} \big)
        \Big] 
        = \frac{1}{n^k} \E_{x_n}\Big[ 
            \sum_{\substack{\mathbf{v} \in U_n^k\\ v_i\ne v_j}} 
            \phi\big( \tfrac{d(\mathbf{v})}{n}, X_{\mathbf{v}} \big)
        \Big] + \widetilde{B}_{x_n,n}\\
        &= 
        \sum_{\sigma \in \mathscr{P}_k} 
        \widetilde{\mathrm{M}}^k_{{x_n},n}[\phi_\sigma \circ D]
        + \widetilde{B}_{x_n,n}.
    \end{align*}
    The error term $\widetilde{B}_{x_n,k}$ is obtained by removing all
    $k$-tuples in the sum which contain duplicate elements of $U_n$. By
    bounding the number of such $k$-tuples by $\widetilde{C}_k
    \abs{U_n}^{k-1}$ (for some combinatorial constant $\widetilde{C}_k$),
    \[
        \widetilde{B}_{{x_n},n} \le \widetilde{C}_k 
        \frac{\norm{\phi}_\infty}{n^k} \E_{x_n}\big[ \abs{U_n}^{k-1} \big].
    \]
    An induction similar to above will prove that $\widetilde{B}_{{x_n}, n} \to
    0$, and thus Proposition~\ref{prop:momentComputation} shows that
    \[
        n \E_{x_n}\big[ \Phi(\bar{\mathcal{U}}_n) \big] 
        \goesto{n \to \infty}
        h(x) \Big( \frac{\Sigma^2}{2} \Big)^{k-1}
        \sum_{\sigma \in \mathscr{P}_k}
        \int 
        \E\big[ \phi_\sigma\big(D(\theta), (X_i)_i \big) \big]
        \widetilde{\Lambda}_k(\diff \theta).
    \]
    The right-hand side is, up to a factor $h(x)\Sigma^2/2$,
    the moment of the marked Brownian CPP with variance $\Sigma^2$
    computed in Proposition~\ref{prop:momentCPP}. We conclude
    again using the method of moments, namely Proposition~\ref{prop:methodMoments}.
    Note that, since the metric space $\bar{\mathcal{U}}_n$ has diameter
    at most $2$, the Gromov-vague and Gromov-weak convergence
    coincide in that case.

    Let us finally argue that these limits hold uniformly. Starting
    from a sequence $(x_n)_{n \ge 1}$, we have shown that, given any
    subsequence of $(x_n)_{n \ge 1}$, we can extract a further
    subsequence such that both limits in \eqref{eq:vagueLimits} hold. This
    shows that, for any metric $d$ that induces the vague convergence,
    \[
        d\big( n \mathscr{L}_x( \bar{\mathcal{T}}_n ),
        h(x_n) \mathscr{L}( \mathcal{T}_{\mathrm{b}, \pi} ) \big)
        \goesto{n\to \infty} 0,
        \qquad
        d\big( n \mathscr{L}_x( \bar{\mathcal{U}}_n ),
        h(x_n) \mathscr{L}( \mathcal{U}_{\mathrm{b}, \pi} ) \big)
        \goesto{n\to \infty} 0.
    \]
    Since $(x_n)_{n \ge 1}$ is arbitrary, this proves that these limits
    are uniform in $x \in E$.
\end{proof}

\begin{proof}[Proof of Theorem~\ref{thm:GromovVague}]
    Point~(i) of the statement has already been proved in
    Proposition~\ref{prop:firstStep}. For Point~(ii),
    Proposition~\ref{prop:firstStep} shows that 
    \[
        \forall x \in E,\quad 
        \frac{n\Sigma^2}{2h(x)} \mathscr{L}_x( \bar{\mathcal{U}}_n )
        \goesto{n \to \infty}
        \mathscr{L}(\mathcal{U}_{\mathrm{b},\pi}),
    \]
    vaguely in the Gromov-weak topology. The result follows from
    \cite[Proposition~3.4]{foutel2024vague}.
\end{proof}

\subsection{Tightness results}
\label{SS:tightness}

This section provides technical lemmas that will allow us to deduce tightness of
the measures. We start by reinforcing the Gromov-vague convergence to a
Gromov-weak one using the estimate on the survival probability. It is
useful for the proofs to recall the notation $Z_m$ for the size of the
population at generation $m \ge 0$.

\begin{corollary} \label{cor:gromovWeak}
    Suppose that the Assumption~\ref{ass:criticality} and
    Assumption~\ref{ass:kolmogorov} hold. Then
    \[
        \mathscr{L}_x( \bar{\mathcal{T}}_n ) 
        \goesto{n \to \infty}
        h(x) \mathscr{L}(\mathcal{T}_{\mathrm{b},\pi})
    \]
    vaguely in the Gromov-weak topology and uniformly for $x \in E$,
    where $\mathcal{T}_{\mathrm{b},\pi}$ is a free Brownian CRT with
    variance $\Sigma^2$ and independent marks $\pi$.
\end{corollary}

\begin{proof}
    Since the convergence holds in the Gromov-vague topology by
    Proposition~\ref{prop:firstStep}, Lemma~\ref{lem:Gv2Gw} shows that it
    is sufficient to prove that, for all $\epsilon > 0$, there is some $R >
    0$ such that
    \[
        \limsup_{n \to \infty} \sup_{x \in E} n \P_x\Big( \sum_{m \ge nR} Z_m \ge \epsilon n^2 \Big) \le \epsilon,
    \]
    Using the survival probability estimate \eqref{eq:kolmogorov}, we
    directly obtain that
    \[
        n \P_x\Big( \sum_{m \ge nR} Z_m \ge \epsilon n^2 \Big) 
        \le n \P_x\big( Z_{nR} > 0 \big) 
        \goesto{n \to \infty}
        \frac{2h(x)}{\Sigma^2 R}
    \]
    uniformly for $x \in E$, hence our result.
\end{proof}

We now prove tightness in the Gromov--Hausdorff-weak topology. Let us
give a verbal description of the ideas before moving to the formal proof.
The tightness criterion that we use (adapted from \cite{athreya16})
relies on controlling the lower mass function $m_{\delta}$ defined in
\eqref{eq:lowerMassFunc}. In words, our task is to show that ``long''
subtrees of $T$ (with height larger then $\delta n$) have
non-vanishing mass (larger than some $\epsilon n^2$). We start by proving
that this property holds for the entire tree $T$, and then use a first
moment argument and the Markov property to show that it also holds for
all subtrees of $T$.

Following the notation of Section~\ref{SS:GromovTopologies}, we denote by
$\tau^{(R)}$ the restriction of a tree $\tau$ to its first $R$
generations and by $\abs{\tau}$ the number of its vertices. 
It will be convenient to introduce a notion of (essential) \emph{height}
of an mmm-space as
\begin{equation} \label{eq:heightDef}
    \forall \mathcal{X} \in \mathscr{X},\quad 
    \h(\mathcal{X}) = \sup \{ d(\rho, x) : (x, e) \in \supp \mu \}.
\end{equation}
We start with a simple lemma.

\begin{lemma} \label{lem:lscHeight}
    The map $\mathcal{X} \mapsto \h(\mathcal{X})$ is lower
    semi-continuous for the Gromov-weak topology.
\end{lemma}

\begin{proof}
    Let $(\mathcal{X}_n)_{n \ge 1}$ be a sequence of mmm-spaces
    converging to $\mathcal{X}$ in the Gromov-weak topology and such that
    $\h(\mathcal{X}_n) \le t$. Then, by Portmanteau's theorem
    \[
        \int_{X \times E} \indic_{\{ d(\rho, x) > t \}} \mu(\diff x, \diff e)
        \le 
        \liminf_{n \to \infty} \int_{X_n \times E} \indic_{\{ d_n(\rho_n, x) > t \}} 
        \mu_n(\diff x, \diff e)
        = 0.
    \]
    Hence $\{ \mathcal{X} : \h(\mathcal{X}) \le t\}$ is closed.
\end{proof}

\begin{lemma} \label{lem:compHeightMass}
    For any $t > 0$, 
    \[
        \adjustlimits
        \lim_{\eta \to 0} \limsup_{n \to \infty} 
        \sup_{x \in E} \P_x( Z_{tn} > 0, \abs{T^{(tn)}} < \eta n^2 ) = 0.
    \]
\end{lemma}

\begin{proof}
    We decompose
    \[
        \P_x( Z_{tn} > 0, \abs{T^{(tn)}} < \eta n^2 ) 
        = 
        \P_x( Z_{tn} > 0) - 
        \P_x( Z_{tn} > 0, \abs{T^{(tn)}} \ge \eta n^2 ).
    \]
    Since $\mathcal{X} \mapsto \h(\mathcal{X})$ is lower semi-continuous
    by Lemma~\ref{lem:lscHeight}, the set $\{ \mathcal{X} :
    \h(\mathcal{X}) > t \}$ is open set for Gromov-weak topology. Therefore, Portmanteau's
    theorem, the vague convergence of $(\bar{\mathcal{T}}_n)_{n \ge 1}$,
    and Kolmogorov's estimate \eqref{eq:kolmogorov} allow us to deduce that
    \begin{multline*}
        \limsup_{n \to \infty}
        \Big( n\P_x( Z_{tn} > 0) - n\P_x( Z_{tn} > 0, \abs{T^{(tn)}} \ge \eta n^2 ) \Big) \\
        \le \frac{2h(x)}{\Sigma^2 t} - h(x) \P( \abs{T^{(t)}_{\mathrm{b}}} \ge \eta,
        \h(\mathcal{T}_{\mathrm{b}}) > t ),
    \end{multline*}
    uniformly for $x \in E$ and where $\mathcal{T}_{\mathrm{b}}$ is a
    free Brownian CRT with variance $\Sigma^2$. By our definition of the
    Brownian CRT, 
    \[
        \lim_{\eta \to 0} \P( \abs{T^{(t)}_\mathrm{b}} \ge \eta,
        \h(\mathcal{T}_\mathrm{b}) > t )
        =
        \P( \h(\mathcal{T}_\mathrm{b}) > t )
        = \frac{2}{\Sigma^2 t},
    \]
    hence the result. (The $\eta \to 0$ limit is obtained by noting that
    the amount of time that a Brownian excursion hitting level $t$ spends
    below level $t$ is positive with probability one.)
\end{proof}

\begin{lemma} \label{lem:tightnessGHP}
    For any $\delta, \eta' > 0$, 
    \[
        \adjustlimits
        \lim_{\eta \to 0} 
        \limsup_{n \to \infty} \sup_{x \in E} 
        n\P_x\big( 
            m_{\delta}(\bar{\mathcal{T}}_n) < \eta, 
            \abs{\bar{\mathcal{T}}_n} > \eta'
        \big)
        = 0.
    \]
\end{lemma}

\begin{proof}
    First, if $\eta < \eta'$, 
    \[
    \{ \abs{\bar{\mathcal{T}}_n} > \eta', m_{\delta}(\bar{\mathcal{T}}_n) < \eta \}
        = 
        \{ \abs{\bar{\mathcal{T}}_n} > \eta', m_{\delta}(\bar{\mathcal{T}}_n) < \eta,
        \h(\bar{\mathcal{T}}_n) > \delta / 2 \},
    \]
    since $\h(\bar{\mathcal{T}}_n) \le \delta / 2$ implies that the diameter of
    the metric space is smaller than $\delta$, in which case
    $m_{\delta}(\bar{\mathcal{T}}_n) = \abs{\bar{\mathcal{T}}_n}$.
    Moreover, since by \eqref{eq:kolmogorov}
    \[
        \lim_{n \to \infty} n \P_x( \h(\bar{\mathcal{T}}_n) > t ) = \frac{2h(x)}{\Sigma^2 t}
    \]
    uniformly for $x \in E$, it is sufficient to show that, for any $K
    \ge 1$ and $\delta > 0$,
    \[
        \adjustlimits \lim_{\eta \to 0} \limsup_{n \to \infty} \sup_{x \in E} 
        n\P_x\big( 
            m_{4\delta}(\bar{\mathcal{T}}_n) < \eta, 
            \h(\bar{\mathcal{T}}_n) 
            \in (\delta, K\delta)
        \big)
        = 0.
    \]

    Now, we construct a family of trees such that (1) each ball of $(T,
    d_T)$ of radius $4\delta n$ contains one such tree and (2) each tree 
    survives for at least $\delta n$ generations. To do so, define 
    for $v \in T$
    \[
        T^{(2\delta n)}_{v}
        \coloneqq 
        \{ w \in T : v \preceq w,\: \abs{w} < \abs{v} + 2\delta n \}.
    \]
    This tree is obtained by cutting the subtree rooted at $v$ at height
    $2\delta n$. For each $u \in T$, let $v_u$ be its ancestor at
    generation $\floor{\delta n} (\floor{ \frac{\abs{u}}{\delta n} } -
    1)$. If $\abs{u} \le \delta n$, we set $v_u = \varnothing$. Clearly,
    \[
        \forall w \in T^{(2\delta n)}_{v_u},\quad 
        d_T(w, u) \le d_T(w, v_u) + d_T(v_u, u) \le 4 \delta n,
    \]
    from which we deduce that
    \[
        T^{(2\delta n)}_{v_u} \subseteq B_{(u,4\delta)}(\bar{\mathcal{T}}_n),
    \]
    the closed ball of $\bar{\mathcal{T}}_n$ with center $u$ and radius $4\delta$.
    The key point is to note that by construction $T^{(2\delta n)}_{v_u}$
    survives at least $n \delta$ generations. If $\abs{u} > \delta n$, $u
    \in T^{(2\delta_n)}_{v_u}$ and $\abs{u} - \abs{v_u} \ge \delta n$. If
    $\abs{u} \le \delta n$, $v_u = \varnothing$ and $T^{(2\delta
    n)}_{v_u} = T^{(2\delta n)}$, which has height at least $n \delta$
    since we restricted ourselves to that event. Therefore,
    \begin{align*}
        \{ m_{4\delta}(\bar{\mathcal{T}}_n) < \eta \}
        &\cap \{ \h(T) \in (\delta n, K \delta n) \} \\
        &= \{ \exists u \in T : \abs{B_{(u,4\delta)}(\bar{\mathcal{T}}_n)} < \eta n^2 \}
        \cap \{ \h(T) \in (\delta n, K \delta n) \} \\
        &\subseteq \{ \exists u \in T : \abs{T^{(2\delta n)}_{v_u}} < \eta n^2 \}
        \cap \{ \h(T) \in (\delta n, K \delta n) \} \\
        &\subseteq
        \bigcup_{k=0}^K \{ \exists v \in T : \abs{T_v^{(2\delta n)}} < \eta n^2,\: 
        \h(T_v^{2\delta n}) > \delta n,\: \abs{v} = k\delta n \}.
    \end{align*}
    By a union bound and the Markov property,
    \begin{align*}
        n\P_x\big( 
            m_{4\delta}(\bar{\mathcal{T}}_n) < \eta, 
            &\h(T) \in (\delta n, K \delta n) 
        \big) \\
        &\le \sum_{k=0}^K n \E_x\Big[ 
            \sum_{\abs{v} = k\delta n} 
            \indic_{\{ \abs{T_v^{(2\delta n)}} < \eta n^2,\: 
                \h(T_v^{(2\delta n)}) > \delta n \}}
            \Big] \\
        &= \sum_{k=0}^K \E_x\Big[ 
            \sum_{\abs{v} = k\delta n} 
            n \P_{X_v}\big( \abs{T^{(2\delta n)}} < \eta n^2,\: 
                Z_{\delta n} > 0 \big)
            \Big] \\
        &\le  \Big(\sup_{y \in E} 
            n \P_y\big( \abs{T^{(2\delta n)}} < \eta n^2,\: Z_{\delta n} > 0 \big) \Big)
            \sum_{k=0}^K \E_x\big[ Z_{k\delta n} \big].
    \end{align*}
    By \eqref{eq:convMean} the second term is bounded uniformly in $x \in
    E$ as $n \to \infty$. Moreover, by Lemma~\ref{lem:compHeightMass}, 
    \begin{equation*}
        \limsup_{n \to \infty} \sup_{y \in E} n \P_y\big( \abs{T^{(2\delta n)}} < \eta n^2,\: Z_{\delta n} > 0 \big) 
        \goesto{\eta \to 0}
        0,
    \end{equation*}
    which leads to our result.
\end{proof}

\subsection{Proof of Theorem~\ref{thm:GHP}}
\label{SS:proofGHP}

\begin{proof}[Proof of Theorem~\ref{thm:GHP}]
    By the Gromov-weak convergence obtained in Corollary~\ref{cor:gromovWeak} 
    and by Lemma~\ref{lem:tightnessGHP}, the assumptions of Lemma~\ref{lem:Gw2GHP} 
    are fulfilled and hence
    \[
        \forall x \in E,\quad 
        n \mathscr{L}_x(\bar{\mathcal{T}}_n) \goesto{n \to \infty}
        h(x) \mathscr{L}(\mathcal{T}_{\mathrm{b}, \pi})
    \]
    vaguely for the Gromov--Hausdorff-weak topology. Note that this
    convergence holds uniformly in $x \in E$. Since each 
    $\bar{\mathcal{T}}_n$ has almost surely full support, this
    convergence also holds vaguely with respect to the
    Gromov--Hausdorff--Prohorov topology, see for instance
    \cite[Remark~5.2]{athreya16}.
\end{proof}

\begin{proof}[Proof of Corollary~\ref{cor:GHPconsequence}]
    We start with point (i). Let $F \colon \mathscr{X} \to \R$ be
    bounded and continuous for the Gromov--Hausdorff-weak topology and
    write
    \begin{equation*}
        n \E_x\big[ F(\bar{\mathcal{T}}_n) \indic_{\{ Z_{tn} > 0 \}} \big] 
        =
        n \E_x\big[ F(\bar{\mathcal{T}}_n) 
        \indic_{\{ Z_{tn} > 0, \abs{\bar{T}^{(t)}_n} > \eta\}} \big]
        +
        n \E_x\big[ F(\bar{\mathcal{T}}_n) 
        \indic_{\{ Z_{tn} > 0, \abs{\bar{T}^{(t)}_n} \le \eta\}} \big].
    \end{equation*}
    By Point~(i) of Theorem~\ref{thm:GHP},
    \begin{align*}
        \lim_{\eta \to 0} \lim_{n \to \infty} 
        n \E_x\big[ F(\bar{\mathcal{T}}_n) 
        \indic_{\{ Z_{tn} > 0, \abs{\bar{T}^{(t)}_n} > \eta\}} \big] 
        &= \lim_{\eta \to 0} 
        h(x) \E\big[ 
            F(\mathcal{T}_{\mathrm{b},\pi}) 
            \indic_{\{ \h(\mathcal{T}_{\mathrm{b},\pi}) > t, \abs{T^{(t)}_{\mathrm{b,\pi}}} > \eta\}} 
        \big]  \\
        &= h(x) \E\big[ 
            F(\mathcal{T}_{\mathrm{b},\pi}) 
            \indic_{\{ \h(\mathcal{T}_{\mathrm{b},\pi}) > t\}} 
        \big] \\
        &= \frac{2h(x)}{\Sigma^2 t} \E\big[ 
            F(\mathcal{T}_{\mathrm{b},\pi})  \:\big|\:
            \h(\mathcal{T}_{\mathrm{b},\pi}) \ge t 
        \big].
    \end{align*}
    By Lemma~\ref{lem:compHeightMass}, the second term can be bounded by
    \begin{align*}
        \adjustlimits \lim_{\eta \to 0} \limsup_{n \to \infty}
        n \E_x\big[ &F(\bar{\mathcal{T}}_n) \indic_{\{ Z_{tn} > 0, \abs{\bar{T}^{(t)}_n} \le \eta\}} \big] \\
        &\le
        \adjustlimits \lim_{\eta \to 0} \limsup_{n \to \infty}
        n\norm{F}_\infty 
        \P_x\big( Z_{tn} > 0, \abs{\bar{T}^{(t)}_n} \le \eta \big) \\
        &= 0.
    \end{align*}
    Combining the two estimates, dividing by the survival probability and
    using \eqref{eq:kolmogorov} leads to the result.

    We move on to point (ii). Let $A$ be a continuity set for the
    free Brownian CRT ($\P(\mathcal{T}_{\mathrm{b}, \pi} \in A) = 0$)
    such that $A \cap \{ \abs{X} \le \eta \} = \emptyset$ for $\eta$
    small enough. Then,
    \begin{equation} \label{eq:convPP}
        \sum_{i=1}^{z_{0,n}} \P_{x_i}(\bar{\mathcal{T}}_n \in A)
        = \frac{1}{n} \sum_{i=1}^{z_{0,n}} n \P_{x_i}(\bar{\mathcal{T}}_n \in A)
        \goesto{n \to \infty}
        \angle{\nu_0, h} \P( \mathcal{T}_{\mathrm{b}, \pi} \in A )
    \end{equation}
    since the initial particle configuration converges weakly and the
    probabilities converge uniformly by the proof of
    Theorem~\ref{thm:GHP}. Viewing the forest as a point measure,
    standard theorems for superposition of independent point processes,
    see \cite[Theorem~11.2.V]{daley2008}, show that
    \eqref{eq:convPP} is enough to deduce that
    \[
        \sum_{i=1}^{z_{0,n}} \delta_{\bar{\mathcal{T}}_{n,i}}
        \goesto{n \to \infty}
        \sum_{i \ge 1} \delta_{\mathcal{T}_{\mathrm{b},i}}
    \]
    in distribution for the vague topology, where
    $(\mathcal{T}_{\mathrm{b},i})_{i \ge 1}$ are the atoms of a Poisson
    point process with intensity $\angle{\nu_0, h} \mathscr{L}(\mathcal{T}_{\mathrm{b}, \pi})$.
    A simple adaptation of standard arguments would show that, 
    if $(\mathcal{T}_i)_{i \ge 1}$ are mmm-spaces such that $\abs{T_1}
    > \abs{T_2} > \dots$, the map $\sum_i \delta_{\mathcal{T}_i}
    \mapsto (\mathcal{T}_i)_{i \ge 1}$ from the space of locally finite point
    measures on $\mathscr{X}$ to the space of sequences of mmm-spaces
    with non-increasing mass is a bijection. This bijection is continuous
    from the vague topology to the product topology, provided that there
    are no ties in the sequence of masses, which occurs almost surely for
    the Brownian CRT.
\end{proof}

\bookmarksetup{startatroot}
\belowpdfbookmark{Acknowledgements}{Acknowledgements}
\section*{Acknowledgement}

The author would like to thank Emmanuel Schertzer for many discussions on
topics closely related to the present work and for his comments on an
earlier version of the manuscript, Magdalen College Oxford for a senior
Demyship, and acknowledges financial support from the Glasstone
Research Fellowship.

\bookmarksetup{startatroot}
\belowpdfbookmark{References}{References}

\bibliographystyle{plain}
\bibliography{many2few}

\appendix

\section{Appendix}

This appendix contains several elementary results on discrete planar
trees.

\begin{lemma} \label{lem:inductiveConstruction}
    The application
    \[
        \tau \mapsto \big(\abs{w}, \big( S_i(\tau) \big)_{i \le \abs{\pi(\tau)}} \big)
    \]
    is a bijection from $\mathscr{T}_k$ to the set
    \[
        \bigcup_{\pi} \N \times 
        \Big( \mathscr{T}_{\abs{\pi_1}} 
        \times \dots \times 
        \mathscr{T}_{\abs{\pi_{\abs{\pi}}}} \Big),
    \]
    where the union runs over all partitions of $[k]$ whose blocks are
    made of consecutive integers, but the partition $\pi = {\{i\}}$.
\end{lemma}

\begin{proof}
    We construct the inverse of the map. Let $\pi = (\pi_i)_{i \le d}$ be
    a non-trivial partition of $[k]$, $n \ge 0$, and $\tau_{1}, \dots, 
    \tau_{d}$ such that $\tau_{i} \in \mathscr{T}_{\abs{\pi_i}}$ and $c_1
    + \dots + c_d = k$. Let $w = (1, \dots, 1)$, with $\abs{w} = n$. We
    construct the tree
    \[
        \tau = \{ v : v \preceq w \} \cup \bigcup_{i = 1}^d
        \{ (w, i, v) : v \in \tau_i \}.
    \]
    Quite clearly, $w$ is the first branch point of $\tau$, it has degree
    $d_w(\tau) = d$, and $S_i(\tau) = \tau_i$ for $i \le d$. 
\end{proof}

\begin{lemma} \label{lem:cppConstruction}
    The application 
    \[
        \tau \mapsto (\ell_i(\tau))_{i \le k},\; (b_i(\tau))_{i<k}
    \]
    is a bijection from $\mathscr{T}_k$ to $\{ (\bm{\uell}, \mathbf{b})
    \in \N^k \times \N^{k-1} : b_i < \ell_i \wedge \ell_{i+1},\; i < k \}$.
\end{lemma}

\begin{proof}
    Let $(\bm{\uell}, \mathbf{b}) \in \N^k \times \N^{k-1}$ be such that 
    $b_i < \ell_i \wedge \ell_{i+1}$. We construct a tree inductively. At
    the first step, $\tau_1 = \{ v : v \preceq v_1 \}$, where $v_1 =
    (1,\dots,1)$ and $\abs{v_1} = \ell_1$. Suppose that $\tau_i$ has been
    constructed. Let $v$ be the right-most (in the planar order) vertex
    of $\tau_i$ at height $\abs{b_i}$. Define $v_{i+1} = (w, 1, \dots,
    1)$ such that $\abs{v_{i+1}} = \ell_{i+1}$, and set 
    \[
        \tau_{i+1} = \tau_i \cup \{ v : v \preceq v_{i+1} \}.
    \]
    Then, $v_{i+1}$ is the $i+1$-th leave and $v_i \wedge v_{i+1} = w$
    has height $b_i$.
\end{proof}

\begin{lemma} \label{lem:factorial2Moment}
    There exists $C_k$ such that, for any finite planar tree $T$, 
    \[
        \Card \{ \mathbf{v} = (v_1,\dots,v_k) \in T^k : \text{$\tau_\mathbf{v}$ has fewer
        than $k$ leaves} \} \le C_k \abs{T}^{k-1} \h(T).
    \]
\end{lemma}

\begin{proof}
    Let $D_k(T)$ denote the quantity that needs to be bounded. The tree 
    $\tau_{\mathbf{v}}$ has strictly less that $k$ leaves if either 
    $\tau_{v_1, \dots, v_{k-1}}$ has strictly less than $k-1$ leaves, or
    if $\tau_{v_1,\dots, v_{k-1}}$ has $k-1$ leaves and $v_k \in 
    \tau_{v_1, \dots, v_{k-1}}$. In the fist case, $v_k$ can be any of
    the $\abs{T}$ vertices of $T$. In the second case, $v_k$ needs to be
    one of the vertices of $\tau_{v_1,\dots, v_{k-1}}$. There are at
    most $(k-1)\h(T)$ such vertices and $\abs{T}^{k-1}$ possible choices
    for $(v_1, \dots, v_{k-1})$. Together, this shows that
    \[
        D_k(T) \le \abs{T} D_{k-1}(T) + \abs{T}^{k-1} (k-1) \h(T).
    \]
    A simple induction now proves the result, given that $D_1(T) = 0$.
\end{proof}

\end{document}